\documentclass[11pt,amsfonts]{article}
\usepackage{graphicx}
\usepackage{latexsym}
\usepackage{amssymb}
\usepackage{amsmath}
\usepackage{enumerate}
\usepackage{color}
\usepackage{bbm}
\usepackage{multirow}
\usepackage{multicol}
\usepackage{caption}
\usepackage{graphicx}

\usepackage{layout}
\newtheorem{prop}{Proposition}
\newtheorem{lemma}{Lemma}
\newtheorem{definition}{Definition}

\newtheorem{theorem}{Theorem}
\newtheorem{remark}{Remark}

\def\qed{ \hfill \vrule width.25cm height.25cm depth0cm\smallskip}
\newcommand{\R}{\mathbb{R}}
\newcommand{\N}{\mathbb{N}}
\newcommand{\Z}{\mathbb{Z}}

\begin{document}
	
	\renewcommand{\thefootnote}{\fnsymbol{footnote}}
	
	\renewcommand{\thefootnote}{\fnsymbol{footnote}}
	\title{Modified weighted power variations of the Hermite process and applications to integrated volatility}
	\author{ A. Ayache\footnote{Université de Lille, 
CNRS, UMR 8524 - Laboratoire Paul Painlev\'e, F-59000 Lille, France. antoine.ayache@univ-lille.fr}, L. Loosveldt\footnote{ \underline{\textbf{Corresponding author:}} Université de Liège, Département de Mathématique -- zone Polytech 1, 12 allée de la Découverte, B\^at. B37, B-4000 Liège. l.loosveldt@uliege.be} and C. A. Tudor\footnote{Université de Lille, 
CNRS, UMR 8524 - Laboratoire Paul Painlev\'e, F-59000 Lille and Bucharest University for Economic Studies, Bucharest, 
Romania. ciprian.tudor@univ-lille.fr}}
	
	\maketitle

\begin{abstract}
We study the asymptotic behaviour of modified weighted power variations of the Hermite process of arbitrary order. By selecting suitable “good” increments and exploiting their decomposition into dominant independent components, we establish a central limit theorem for weighted 
	$p$-variations using tools from Stein–Malliavin calculus. Our results extend previous works on modified quadratic and wavelet-based variations to general powers and to weighted settings, with explicit bounds in Wasserstein distance. We further apply these limit theorems to construct asymptotically Gaussian estimators of integrated volatility in Hermite-driven models, thereby extending fBm-based methods to non-Gaussian settings. The last part of our work contains numerical simulations which  illustrate  the practical performance of the proposed estimators.
\end{abstract}


	
\section{Introduction}
The study of Hermite and related processes has developed intensively in the recent decades. These stochastic processes enjoy nice properties: they are self-similar, have stationary increments and exhibit long-range dependence. These characteristics give them an important potential for practical applications. The class of Hermite processes includes the fractional Brownian motion which is the only Gaussian process in this class. We refer to the monographs \cite{PiTa-book}  and \cite{T3} for a detailed exposition on Hermite and related processes. 

One of the important issues in the analysis of Hermite processes is the statistical estimation of its self-similarity index (or Hurst parameter), which characterizes  many of its  properties. A series of papers (\cite{CTV}, \cite{TV}, \cite{BaTu}, \cite{CRTT}) has proposed estimators for the Hurst index of the Hermite processes, by using power variations or wavelet-based methods.  A common fact for all these estimators is that, although they are consistent, they are not asymptotically Gaussian and this constitutes a serious drawback for statistical estimation. A similar phenomenon can be noticed for other parameter estimators in Hermite-driven models (see e.g. \cite{NT}, \cite{Tran}).

To overcome the asymptotic non-normality, a new strategy has been employed in \cite{AT}, by following an idea from \cite{A}. Namely, the reference \cite{AT} has defined a so-called modified quadratic variation of the Hermite process (of an arbitrary order), by using only some well-chosen increments of it. Each of these special increments can be split into a dominant part and a negligible one, moreover the  dominant parts constitute a sequence of independent random variables. Thanks to these nice properties, \cite{AT} has managed to show that the modified quadratic variation satisfies a Central Limit Theorem, and then has constructed from it a new strongly consistent estimator for the Hurst parameter of the Hermite processes which, compared to other estimators, offers the significant advantage to be asymptotically Gaussian. A related asymptotically Gaussian estimator based on the modified wavelet variations of the Hermite process has been defined and analyzed in \cite{LT}. 

In this work, our purpose is to generalize the results obtained in \cite{AT} and \cite{LT} by analyzing the limit behavior of the weighted modified $p$-variations of the Hermite process $(Z_{t}, t\geq 0)$  for  any integer order $p \geq 2$, i.e. the sequence $(U_{N,p}(Z), N\geq 1)$ given by (\ref{unpz}). The weight $h$ in (\ref{unpz}) is a deterministic function with suitable properties.  Via a similar idea of choosing ``good'' increments of the Hermite process and by using the techniques of the Stein-Malliavin calculus, we prove that the weighted modified power variation of the Hermite process converges to a constant times the Wiener integral of the weight function $h$ with respect to a Brownian motion independent of the underlying Hermite process. We start by considering the case of the  (non-weighted) modified power variation (i.e. $h$ is identically equal to one) and then we extend the result to the weighted case by using, among others, some new findings  concerning  the Wasserstein  distance. 

We apply our result to estimate the integrated volatility in a simple model with Hermite noise. The estimation of the integrated volatility via power variations  constitutes a topic of interest in financial mathematics. We mention the articles \cite{BCP}, \cite{Betal}, \cite{PV} or \cite{Ja1} for 
semimartingale-type models and the papers \cite{CNW} or \cite{BCP} for models driven by the fractional Brownian motion. Our work extends the results obtained in the fBm case to non-Gaussian settings.

We organize our paper as follows. Section 2 contains preliminaries on Hermite processes, including their representation, basic properties, and the definition of the special increments that play a central role in our analysis. This section also contains some establishes key independence and moment estimates for these increments. In Section 3 we include  the proof of some properties of strongly independent random variables relying on tools from  Wiener chaos. We then investigate, in Section 4, the modified (non-weighted) power variations, proving a law of large numbers and a central limit theorem through a combination of Stein–Malliavin arguments and independence properties. The subsequent Section 5 extends these results to the weighted setting, deriving quantitative Wasserstein bounds and introducing a general limit theorem for weighted modified power variations. We also apply the theoretical results to integrated volatility estimation in Hermite-driven models and present numerical simulations that validate the theoretical findings and illustrate the performance of the proposed estimators. Section 6 is the Appendix where we present the basic elements of the Malliavin calculus.

Throughout, we denote by $C_a$ a strictly positive constant that depends on the parameter $a$.
\section{Preliminaries}
In this section we present various facts related to  the Hermite process. The first part contains its  definition and basic properties and in the second part we describe some special increments of this stochastic process, which play a crucial role in our construction.

\subsection{The Hermite process and its special increments}\label{sec21}
We denote by $ (Z_{t}, t\geq 0)$ the Hermite process of order $q\geq 1$ and with self-similarity index $ H\in \left( \frac{1}{2}, 1\right)$. It can be expressed, for every $t\geq 0$, as 
	\begin{equation}\label{hermite1}
	Z_{t}= d(q,H) \int_{ \mathbb{R} ^{q}} \left( \int_{0} ^{t} f_{u} (y_{1},...,y_{q}) du \right) dB(y_{1})...dB(y_{q}), 
\end{equation}
where $(B(y), y \in \mathbb{R})$ is a Wiener process and, for every $y_{1},..., y_{q}\in \mathbb{R}$, 
\begin{equation}
	\label{fu}
	f_{u}(y_{1},...,y_{q}) =d(q,H) (u-y_{1})_{+} ^{-\left(\frac{1}{2}+ \frac{1-H}{q}\right)}\ldots (u-y_{q})_{+} ^{-\left(\frac{1}{2}+ \frac{1-H}{q}\right)},
\end{equation}
and $d(q,H)$ is the strictly positive normalizing constant which ensures that $ \mathbf{E} Z_{t} ^{2}= t ^{2H}$ for every $t\geq 0$. Alternatively, we can write 
\begin{equation*}
	Z_{t}= I_{q}(L_{t}), \hskip0.5cm t\geq 0,
\end{equation*}
where $ I_{q}$ denotes the multiple stochastic integral of order $q$ with respect to $B$ and
\begin{equation*}
	L_{t} (y_{1},..., y_{q})= \int_{0} ^{t} f_{u}(y_{1},..., y_{q}) du,
\end{equation*}
for every $y_{1},..., y_{q} \in \mathbb{R}$.  The stochastic process $ (Z_{t}, t\geq 0)$ is $H$-self-similar and it has stationary increments. Consequently, its covariance function reads 
\begin{equation*}
	\mathbf{E}  Z_{t} Z_{s}= \frac{1}{2} \left( t ^{2H}+ s^{2H}-\vert t-s\vert ^{2H}\right), \mbox{ for every } s,t\geq 0. 
\end{equation*}
The increments of $Z$ satisfy, for $s, t\geq 0$ and $ p\geq 1$, 
\begin{equation*}
	\mathbf{E}  \left| Z_{t} -Z_{s} \right| ^{p} = \mathbf{E}  \vert Z_{1}\vert ^{p}\vert t-s\vert ^{Hp}. 
\end{equation*}
Hence, for every fixed $\delta \in (0, H)$, sample paths of the Hermite process are, on each compact subinterval of $\mathbb{R}_+$,
$\delta$-H\"older continuous functions.

In the sequel we will use the Wiener integral with respect to the Hermite process. Let $\vert \mathcal{H}\vert $ be the space of measurable functions $f: \mathbb{R} \to \mathbb{R}$ such that 
\begin{equation*}
	\int_{\mathbb{R}} \int_{\mathbb{R}} \vert f(u)\vert \cdot \vert f(v) \vert \cdot \vert u-v\vert ^{2H-2} du dv < \infty. 
\end{equation*}
If $ f\in \vert \mathcal{H}\vert$, it is possible to define the Wiener integral of $f$ with respect to the Hermite process $Z$ (also called the Wiener-Hermite integral), denoted by
\begin{equation*}
	\int_{\mathbb{R}} f(u) dZ_{u}.
\end{equation*}
We refer to   e.g.  the survey \cite{T3} for the details of this construction. The above object is well-defined as a random variable in $ L ^{2} (\Omega)$ and it satisfies the following isometry property: if $f,g\in \vert \mathcal{H}\vert$, then 

\begin{equation}\label{iso-w}
	\mathbf{E} \Bigg( \left( \int_{\mathbb{R}} f(u) dZ_{u}\right) \left( \int_{\mathbb{R}} g(u) dZ_{u}\right)\Bigg)  =H(2H-1) \int_{\mathbb{R}}\int_{\mathbb{R}}f(u) g(v) \vert u-v\vert ^{2H-2}dudv. 
\end{equation}

\subsection{The special increments of the Hermite process}\label{sec22}
We will deal with some special increments of the Hermite process introduced in \cite{AT}  (see also \cite{LT} for an approach based on wavelet coefficients). They are defined as follows. We first fix a real number $\gamma \in (0, 1)$.  For every integers $N\geq 1$ and $l=1,..., 2 ^{[N ^{\gamma}]}$, we consider the following increments of length $ 2 ^{-N}$ located at the points $\frac{l}{ 2 ^{[ N ^{\gamma}]}}$ of the  unit interval $[0, 1]$:
	\begin{equation}\label{dz}
	\Delta Z_{l, N}= Z _{ \frac{l}{2 ^{ [N ^{\gamma}]}}+2 ^{-N}}-  Z _{ \frac{l}{2 ^{ [N ^{\gamma}]}}}.
\end{equation}
The key observation made in the reference \cite{A} is that each of the above increments can be decomposed in a negligible part with small variance and a dominant part, and their  dominants parts are mutually independent for $ l=1,..., 2 ^{ [N ^{\gamma}]}$. Let us recall this decomposition.  First, by (\ref{hermite1}), we can write 
\begin{eqnarray*}
	\Delta Z_{l, N}&=& \int_{\mathbb{R} ^{q}} \mathbbm{1}_{ \left( -\infty, \frac{l}{2 ^{[N ^{\gamma}]}}+2 ^{-N}\right)^q} \left( \int_{ \frac{ l}{2 ^{[N ^{\gamma}]}}}^{  \frac{ l}{2 ^{[N ^{\gamma}]}}+2 ^{-N}} f_{u} (y_{1},...,y_{q}) du \right) dB(y_{1})...dB(y_{q}).
\end{eqnarray*}
Notice that the indicator function in the previous integral, as well as the indicator functions which will appear in the sequel, depend on $(y_{1},...,y_{q})$.

	Let us now consider a fixed real number  $\beta$ such that $\gamma<\beta <1$. From now on, we always assume that $N\ge 2^{(1-\beta)^{-1}}$. We decompose the above increment $ \Delta Z_{l, N}$ in the following way
\begin{eqnarray}
&&	\Delta Z_{l, N}\nonumber \\&=& \int_{\mathbb{R} ^{q}} \mathbbm{1}_{ \left(\frac{l}{2 ^{[N ^{\gamma}]}}-\frac{ 2 ^{[N ^{\beta}]}}{2 ^{N}}+ 2 ^{-N} , \frac{l}{2 ^{[N ^{\gamma}]}}+2 ^{-N}\right)^q}\left( \int_{ \frac{ l}{2 ^{[N ^{\gamma}]}}}^{  \frac{ l}{2 ^{[N ^{\gamma}]}}+2 ^{-N}} f_{u} (y_{1},...,y_{q}) du \right) dB(y_{1})...dB(y_{q})\nonumber \\
	&&+  \int_{\mathbb{R} ^{q}} \mathbbm{1}_{ \overline{\left(\frac{l}{2 ^{[N ^{\gamma}]}}-\frac{ 2 ^{[N ^{\beta}]}}{2 ^{N}}+ 2 ^{-N} , \frac{l}{2 ^{[N ^{\gamma}]}}+2 ^{-N}\right)^{q}}}\left( \int_{ \frac{ l}{2 ^{[N ^{\gamma}]}}}^{  \frac{ l}{2 ^{[N ^{\gamma}]}}+2 ^{-N}} f_{u} (y_{1},...,y_{q}) du \right) dB(y_{1})...dB(y_{q})\nonumber\\
	&=& \tilde{\Delta} Z_{l, N} + \check{\Delta}Z_{l, N},\label{deco}
\end{eqnarray} 
where 
\begin{eqnarray*}
&&	\overline{\left(\frac{l}{2 ^{[N ^{\gamma}]}}-\frac{ 2 ^{[N ^{\beta}]}}{2 ^{N}}+ 2 ^{-N} , \frac{l}{2 ^{[N ^{\gamma}]}}+2 ^{-N}\right)^{q}}\\
&&= \left( -\infty, \frac{l}{2 ^{[N ^{\gamma}]}}+2 ^{-N}\right)^{q}\setminus \left(\frac{l}{2 ^{[N ^{\gamma}]}}-\frac{ 2 ^{[N ^{\beta}]}}{2 ^{N}}+ 2 ^{-N} , \frac{l}{2 ^{[N ^{\gamma}]}}+2 ^{-N}\right)^{q}.
\end{eqnarray*}
We used the notation
\begin{equation}
	\label{dt}
	\tilde{\Delta} Z_{l, N} = \int_{\mathbb{R} ^{q}} \mathbbm{1}_{ \left(\frac{l}{2 ^{[N ^{\gamma}]}}-\frac{ 2 ^{[N ^{\beta}]}}{2 ^{N}}+ 2 ^{-N} , \frac{l}{2 ^{[N ^{\gamma}]}}+2 ^{-N}\right)^q}\left( \int_{ \frac{ l}{2 ^{[N ^{\gamma}]}}}^{  \frac{ l}{2 ^{[N ^{\gamma}]}}+2 ^{-N}} f_{u} (y_{1},...,y_{q}) du \right) dB(y_{1})...dB(y_{q})
\end{equation}
and
\begin{equation}
	\label{dc}
	\check{\Delta}Z_{l, N}= \int_{\mathbb{R} ^{q}} \mathbbm{1}_{ \overline{\left(\frac{l}{2 ^{[N ^{\gamma}]}}-\frac{ 2 ^{[N ^{\beta}]}}{2 ^{N}}+ 2 ^{-N} , \frac{l}{2 ^{[N ^{\gamma}]}}+2 ^{-N}\right)^{q}}}\left( \int_{ \frac{ l}{2 ^{[N ^{\gamma}]}}}^{  \frac{ l}{2 ^{[N ^{\gamma}]}}+2 ^{-N}} f_{u} (y_{1},...,y_{q}) du \right) dB(y_{1})...dB(y_{q}).
\end{equation}

Actually, the assumption that the integer $N$ is large enough so that $N\ge 2^{(1-\beta)^{-1}}$  ensures that the intervals $\left(\frac{l}{2 ^{[N ^{\gamma}]}}-\frac{ 2 ^{[N ^{\beta}]}}{2 ^{N}}+ 2 ^{-N} , \frac{l}{2 ^{[N ^{\gamma}]}}+2 ^{-N}\right)$, $l=1,..., 2 ^{[N ^{\gamma}]}$ are disjoint. 

\begin{remark}
Let us emphasise that we do not use exactly the same increments \eqref{dz}, nor their decomposition \eqref{deco}, as in the article \cite{AT}. This is because, in \cite{AT}, the estimator was constructed from increments of length $2^{-N}$ located at the points $\frac{l\,[ 2^{N^\beta} ]}{2^N}$ for $l \in \mathbb{N}\cap\left [ 1, \frac{2^N}{[2^{N^\beta} ]} \right ] \cap [1, [ 2^{N^\gamma} ] ]$ and $0<\gamma<\beta<1$. Proceeding in this way, the intervals selected for the variation were not evenly distributed over the whole interval $[0,1]$, since only the first $[ 2^{N^\gamma} ]$ were chosen. This had no effect on the estimation of the Hurst index. In the present work, however, we found this inconvenient for handling the introduction of a weight in the variations, notably for the proof of Theorem \ref{tt3bi}. For this reason, we have defined a new method for selecting increments.
\end{remark}

We next derive some important estimates concerning the moments of the quantities $ \tilde{\Delta}Z _{l, N} $ and $\check{\Delta}Z_{l, N}$ given by (\ref{dt}) and (\ref{dc}), respectively. The first result has been essentially proven in Lemma 2.3 in \cite{A}. Since we slightly changed the definition of the random variables $\tilde{\Delta}Z _{l, N} $ and $\check{\Delta}Z_{l, N}$, we choose to include the proof below.

\begin{lemma}\label{ll1}
	For $ N\geq 1$ and $ l=1,2,..., 2 ^{[N ^{\gamma}]}$, let $	\check{\Delta}Z_{l, N}$ be given by (\ref{dc}).  Then
	\begin{equation*}
		\mathbf{E}  	\check{\Delta}Z_{l, N} ^{2} \leq C(H,q) 2 ^{-2HN} 2 ^{[N ^{\beta} ]\frac{2H-2}{q}},
	\end{equation*}
	for every $ N\geq 1$ and $ l=1,2,..., 2 ^{[N ^{\gamma}]}$. 
\end{lemma}
{\bf Proof: } Since the integrand in (\ref{dc}) is a symmetric function in $(y_{1},...,y_{q})$, by using the isometry  property of the multiple stochastic integrals (\ref{iso}),  we get
\begin{eqnarray*}
	\mathbf{E}  	\check{\Delta}Z_{l, N} ^{2} &=& q! \int_{\mathbb{R} ^{q}} dy_{1}...dy_{q} \mathbbm{1}_{ \overline{\left(\frac{l}{2 ^{[N ^{\gamma}]}}-\frac{ 2 ^{[N ^{\beta}]}}{2 ^{N}}+ 2 ^{-N} , \frac{l}{2 ^{[N ^{\gamma}]}}+2 ^{-N}\right)^{q}}}\\
	&&\int_{ \frac{ l}{2 ^{[N ^{\gamma}]}}}^{  \frac{ l}{2 ^{[N ^{\gamma}]}}+2 ^{-N}} f_{u} (y_{1},...,y_{q}) du \int_{ \frac{ l}{2 ^{[N ^{\gamma}]}}}^{  \frac{ l}{2 ^{[N ^{\gamma}]}}+2 ^{-N}} f_{v} (y_{1},...,y_{q}) dv\\
	&\leq & q! \cdot q\int_{- \infty} ^ { \frac{l}{2 ^{[N ^{\gamma}]}}-\frac{ 2 ^{[N ^{\beta}]}}{2 ^{N}}+ 2 ^{-N}}dy_{1} \int_{\mathbb{R} ^{q-1}} dy_{2}...dy_{q} \\
	&&\int_{ \frac{ l}{2 ^{[N ^{\gamma}]}}}^{  \frac{ l}{2 ^{[N ^{\gamma}]}}+2 ^{-N}}\int_{ \frac{ l}{2 ^{[N ^{\gamma}]}}}^{  \frac{ l}{2 ^{[N ^{\gamma}]}}+2 ^{-N}}dudv f_{u} (y_{1},...,y_{q}) f_{v} (y_{1},...,y_{q}).
\end{eqnarray*}
Now, by using the identity, for all real numbers $-1<a<-\frac{1}{2}$ and $u\ne v$,
\begin{equation*}
	\int_{\mathbb{R}} (u-y)_{+} ^{a} (v-y)_{+} ^{a} dy= \beta (-1-2a, a+1) \vert u-v\vert ^{2a+1},
\end{equation*}
where $\beta$ stands for the beta function given by $\beta (a, b)= \int_{0} ^{1} r ^{a-1} (1-r) ^{b-1}dr$ for $a,b>0$, we get 
\begin{eqnarray*}
	&& \int_{\mathbb{R} ^{q-1}} dy_{2}...dy_{q}  f_{u}(y_{1},..., y_{q})  f_{v}(y_{1},..., y_{q}) \\
	&=& C(q,H)  (u-y_{1})_{+} ^{-\left(\frac{1}{2}+ \frac{1-H}{q}\right)}(v-y_{1})_{+}^{-\left(\frac{1}{2}+ \frac{1-H}{q}\right)}\vert u-v\vert ^{(2H-2) \frac{q-1}{q}}.
\end{eqnarray*}
Therefore,
\begin{eqnarray*}
	\mathbf{E}  	\check{\Delta}Z_{l, N} ^{2} &\leq &C(q,H)\int_{ \frac{ l}{2 ^{[N ^{\gamma}]}}}^{  \frac{ l}{2 ^{[N ^{\gamma}]}}+2 ^{-N}}\int_{ \frac{ l}{2 ^{[N ^{\gamma}]}}}^{  \frac{ l}{2 ^{[N ^{\gamma}]}}+2 ^{-N}}dudv \int_{- \infty} ^ { \frac{l}{2 ^{[N ^{\gamma}]}}-\frac{ 2 ^{[N ^{\beta}]}}{2 ^{N}}+ 2 ^{-N}}dy_{1} \\
	&& (u-y_{1})_{+} ^{-\left(\frac{1}{2}+ \frac{1-H}{q}\right)}(v-y_{1})_{+}^{-\left(\frac{1}{2}+ \frac{1-H}{q}\right)}\vert u-v\vert ^{(2H-2) \frac{q-1}{q}}\\
	&\leq & C(q,H) \int_{ \frac{ l}{2 ^{[N ^{\gamma}]}}}^{  \frac{ l}{2 ^{[N ^{\gamma}]}}+2 ^{-N}}\int_{ \frac{ l}{2 ^{[N ^{\gamma}]}}}^{  \frac{ l}{2 ^{[N ^{\gamma}]}}+2 ^{-N}}dudv \int_{- \infty} ^ { \frac{l}{2 ^{[N ^{\gamma}]}}-\frac{ 2 ^{[N ^{\beta}]}}{2 ^{N}}+ 2 ^{-N}}dy_{1} \\
	&&\left( \frac{l} { 2 ^{[N ^{\gamma}]}}-y_{1}\right) ^{-\left(\frac{1}{2}+ \frac{1-H}{q}\right)}\left( \frac{l} { 2 ^{[N ^{\gamma}]}}-y_{1}\right) ^{-\left(\frac{1}{2}+ \frac{1-H}{q}\right)}\vert u-v\vert ^{(2H-2) \frac{q-1}{q}}\\
	&\leq & C(q,H)\left( \frac{ 2 ^{[N ^{\beta}]}}{2 ^{N}}\right) ^{\frac{ 2H-2}{q}}  \int_{ \frac{ l}{2 ^{[N ^{\gamma}]}}}^{  \frac{ l}{2 ^{[N ^{\gamma}]}}+2 ^{-N}}\int_{ \frac{ l}{2 ^{[N ^{\gamma}]}}}^{  \frac{ l}{2 ^{[N ^{\gamma}]}}+2 ^{-N}}dudv \vert u-v\vert ^{(2H-2) \frac{q-1}{q}}\\
	&=&C(q,H)\left( \frac{ 2 ^{[N ^{\beta}]}}{2 ^{N}}\right) ^{\frac{ 2H-2}{q}} \left( \frac{1}{ 2 ^{N}}\right) ^{(2H-2)\frac{q-1}{q}+2}= C(q,H) 2 ^{-2HN}  2 ^{[N ^{\beta}] \frac{ 2H-2}{q}}.
\end{eqnarray*}
\qed
 
Concerning the moments of the part denoted by $\tilde{\Delta}Z_{l, N}$, we have the following result. The case $p=2$ has been treated in \cite{AT}.

\begin{lemma}\label{ll2}
	For $N\geq 1$ and $l=1,..., 2 ^{[N ^{\gamma}]}$, let $ \tilde{\Delta}Z _{l, N}$ be given by (\ref{dt}). Then, for every integers $p\geq 2$ and $l=1,..., 2 ^{[N ^{\gamma}]}$,
	\begin{equation}\label{15i-1}
		\left| 2 ^{pHN} \mathbf{E}  (\tilde{\Delta}Z _{l, N})^{p} - \mathbf{E}  Z_{1} ^{p}\right| \leq C_p 2 ^{[N ^{\beta}] \frac{ H-1}{q}}.
	\end{equation}
	In particular, 
	\begin{equation*}
		2 ^{pHN} \mathbf{E}  ( \tilde{\Delta} Z_{l, N} )^{p} \to _{N\to \infty} \mathbf{E}  Z_{1} ^{p},
	\end{equation*}
	for every $l=1,..., 2 ^{[N ^{\gamma}]}$.
\end{lemma}
{\bf Proof: } Fix $p\geq 2$ and $l=1,..., 2 ^{[N ^{\gamma}]}$. By the self-similarity and the stationarity of the increments of the Hermite process, 
\begin{equation*}
	\mathbf{E}  (\Delta Z _{l,N})^{p} = 2^{-pHN}\mathbf{E}  Z_{1} ^{p}.
\end{equation*}
By using (\ref{deco}), we can then write 
\begin{equation*}
	\mathbf{E}  (\tilde{\Delta} Z _{l,N} )^{p}+ \sum _{j=0}^{p-1} \binom{p}{j} 	\mathbf{E}  (\tilde{\Delta} Z _{l,N} )^{j}	(\check{\Delta} Z _{l,N} )^{p-j}= 2^{-pHN}\mathbf{E}  Z_{1} ^{p},
\end{equation*}
so 
\begin{eqnarray*}
	\left| 2 ^{pHN} \mathbf{E}  (\tilde{\Delta} Z _{l,N} )^{p} - \mathbf{E}  Z_{1}^{p}\right|&=& 2 ^{pHN} \left| \sum _{j=0}^{p-1} \binom{p}{j} 	\mathbf{E}  (\tilde{\Delta} Z _{l,N} )^{j}	 (\check{\Delta} Z _{l,N} )^{p-j}\right| \\
	&=& 2 ^{pHN} \left| \sum _{j=0}^{p-2} \binom{p}{j} 	\mathbf{E}  (\tilde{\Delta} Z _{l,N} )^{j}	 (\check{\Delta} Z _{l,N} )^{p-j}\right| ,
\end{eqnarray*}
the summand with $j=p-1$ vanishing by Lemma 2 in \cite{AT}. Hence, it follows from Cauchy-Schwarz's inequality that
\begin{eqnarray*} 
	\left| 2 ^{pHN} \mathbf{E}  (\tilde{\Delta} Z _{l,N} )^{p} - \mathbf{E}  Z_{1}^{p}\right|	&\leq & 2 ^{pHN} \sum _{j=0}^{p-2} \binom{p}{j} \left( \mathbf{E}  (\tilde{\Delta} Z _{l,N} )^{2j}\right) ^{\frac{1}{2}}\left( \mathbf{E}  (\check{\Delta} Z _{l,N} )^{2(p-j)}\right) ^{\frac{1}{2}}\\
	&\leq & C_{p} 2 ^{pHN} \sum _{j=0}^{p-2}\left( \mathbf{E}  (\tilde{\Delta} Z _{l,N} )^{2}\right) ^{\frac{j}{2}}\left( \mathbf{E}  (\check{\Delta} Z _{l,N} )^{2}\right) ^{\frac{p-j}{2}},
\end{eqnarray*}
where the last bound is obtained via the hypercontractivity property (\ref{hyper}) of multiple stochastic integral. Finally, one can derive from Lemma \ref{ll1} and the inequality $\mathbf{E}  (\tilde{\Delta} Z _{l,N} )^{2}\le \mathbf{E}  (\Delta Z _{l,N} )^{2}=2^{-2HN}$ that
\begin{eqnarray*}
	\left| 2 ^{pHN} \mathbf{E}  (\tilde{\Delta} Z _{l,N} )^{p} - \mathbf{E}  Z_{1}^{p}\right|&\leq &C_{p} 2 ^{pHN} \sum _{j=0}^{p-1}  2 ^{-jHN} 2 ^{-(p-j)HN} \left( 2 ^{ [N ^{\beta}] \frac{2H-2}{q}}\right) ^{\frac{p-j}{2}}\\
	&\leq & C_{p} 2 ^{[N ^{\beta}] \frac{H-1}{q}}.
\end{eqnarray*}
\qed

The following properties of the random variables $\tilde{\Delta}Z_{l, N}$ and $\check{\Delta}Z_{l,N}$ play an important role in the sequel. 

\begin{prop}\label{pp1}
For $N\geq 2^{(1-\beta)^{-1}}$ and $l=1,..., 2 ^{[N ^{\gamma}]}$, let $\tilde{\Delta}Z_{l, N}$ and $\check{\Delta}Z_{l,N}$  be defined by (\ref{dt}) and (\ref{dc}), respectively. Then 
	\begin{enumerate}
		\item The random variables $\tilde{\Delta}Z_{l, N}, l=1,..., 2 ^{[N ^{\gamma}]}$ are independent and identically distributed. 
		
		\item  The random variables $\check{\Delta}Z_{l, N}, l=1,..., 2 ^{[N ^{\gamma}]}$ are identically distributed. 
	\end{enumerate}
\end{prop}
{\bf Proof: } To prove point 1., we recall the main result in \cite{UZ}: the multiple stochastic integrals $I_{n_{i}}(f_{i}), i=1,..., d$ (where $f_{i}\in L^{2}(\mathbb{R} ^{n_{i}})$  are symmetric functions) are independent if and only if for every $i,j=1,...,d$ with $i\not=j$, 
\begin{equation*}
	f_{i}\otimes _{1}f_{j}=0 \mbox{ almost everywhere on } \mathbb{R} ^{n_{i}+ n_{j}-2}. 
\end{equation*}
We can express $\tilde{\Delta}Z_{l, N}$ as a multiple stochastic integral of order $q$ in the following way
\begin{equation*}
	\tilde{\Delta}_{l, N}=I_{q} (g_{l,N} \mathbbm{1}_{A_{l, N}} ^{\otimes q}),
\end{equation*}
with
\begin{equation*}
	A_{l, N}= \left(\frac{l}{2 ^{[N ^{\gamma}]}}-\frac{ 2 ^{[N ^{\beta}]}}{2 ^{N}}+ 2 ^{-N} , \frac{l}{2 ^{[N ^{\gamma}]}}+2 ^{-N}\right),
\end{equation*}
and, for every $y_{1},..,y_{q}\in \mathbb{R}$,
\begin{equation*}
	g_{l,N}(y_{1},..., y_{q})=  \int_{ \frac{ l}{2 ^{[N ^{\gamma}]}}}^{  \frac{ l}{2 ^{[N ^{\gamma}]}}+2 ^{-N}} f_{u} (y_{1},...,y_{q}) du,
\end{equation*}
with $f_{u}$ given by (\ref{fu}). To obtain the conclusion of point 1., it suffices to show that for every $k,l=1,..., 2 ^{[N ^{\gamma}]}$ with $k\not=l$, 
\begin{equation*}
	\left( g_{l,N} \mathbbm{1}_{A_{l, N}} ^{\otimes q}\right) \otimes _{1} \left( g_{k,N}\mathbbm{1}_{A_{k, N}} ^{\otimes q}\right)=0 \mbox{ almost everywhere on } \mathbb{R} ^{ 2q-2}. 
\end{equation*}
By the definition of the contraction, for $y_{1},..., y_{2q-2}\in \mathbb{R}$,
\begin{eqnarray*}
	&&	\left( g_{l,N} \mathbbm{1}_{A_{l, N}} ^{\otimes q}\right) \otimes _{1} \left( g_{k,N} \mathbbm{1}_{A_{k, N}} ^{\otimes q}\right)(y_{1},..., y_{2q-2})\\
	&=& \int_{\mathbb{R}}dx g_{l, N} (y_{1},..., y_{q-1}, x) \mathbbm{1}_{ A_{l, N}}^{\otimes q-1}(y_{1},..., y_{q-1})\mathbbm{1}_{ A_{l, N}}(x) \\
	&& g_{k, N} (y_{q},..., y_{2q-2}, x)\mathbbm{1}_{ A_{l, N}}^{\otimes q-1}(y_{q},..., y_{2q-2})\mathbbm{1}_{ A_{l, N}}(x)
\end{eqnarray*}
and the latter integral vanishes since, for $k\not=l$, $A_{l, N}$ and $ A_{k, N}$ are disjoint sets when $N\geq 2^{(1-\beta)^{-1}}$. 

To show that the random variables $\tilde{\Delta}Z_{l, N}, l=1,..., 2 ^{[N ^{\gamma}]}$ are identically distributed, we write (with the notation $ = ^{(d)}$ for the equality in distribution),
\begin{eqnarray*}
&&	\tilde{\Delta}Z_{l, N}\\
&=&\int_{\mathbb{R} ^{q}} 1_{ \left(\frac{l}{2 ^{[N ^{\gamma}]}}-\frac{ 2 ^{[N ^{\beta}]}}{2 ^{N}}+ 2 ^{-N} , \frac{l}{2 ^{[N ^{\gamma}]}}+2 ^{-N}\right)}^{\otimes q} \left( \int_{ \frac{ l}{2 ^{[N ^{\gamma}]}}}^{  \frac{ l}{2 ^{[N ^{\gamma}]}}+2 ^{-N}} f_{u} (y_{1},...,y_{q}) du \right) dB(y_{1})...dB(y_{q})\\
	&=^{(d)} & \int_{\mathbb{R} ^{q}} dB(y_{1})...dB(y_{q})\mathbbm{1}_{ \left( -\frac{ 2 ^{[N ^{\beta}]}}{2 ^{N}}, 0\right)}^{\otimes q}(y_{1},..., y_{q})\\
	&&\left( \int_{ \frac{ l}{2 ^{[N ^{\gamma}]}}}^{  \frac{ l}{2 ^{[N ^{\gamma}]}}+2 ^{-N}} f_{u} (y_{1}+ \frac{l}{2 ^{[N ^{\gamma}]}}+2^{-N},...,y_{q}+ \frac{l}{2 ^{[N ^{\gamma}]}}+2^{-N}) du \right)\\
	&=&  \int_{\mathbb{R} ^{q}} \mathbbm{1}_{ \left( -\frac{ 2 ^{[N ^{\beta}]}}{2 ^{N}}, 0\right)}^{\otimes q}(y_{1},..., y_{q})\left( \int_{-2 ^{-N}}^{0} f_{u} (y_{1},..., y_{q}) du \right) dB(y_{1})...dB(y_{q}),
\end{eqnarray*}
and the last quantity is independent of $l=1,..., 2 ^{[N ^{\gamma}]}$.  A similar argument can be used to show point 2. \qed

\section{Strong independence}
Let us recall the notion of strong independence for random variables. It will be used intensively in the next section in order to prove our main results. Let $F$ and $G$ be two random variables that can be expanded in $L^2(\Omega)$ as $ F=\sum_{ n\geq 0} I_{n}(f_{n})$ and $ G= \sum_{n\geq 0} I_{n} (g_{n})$ where  $f_{n}, g_{n} \in L ^{2} (\mathbb{R}^{n})$ are symmetric functions for every $n\geq 1$.  We say that $F$ and $G$ are strongly independent if $I_{m}(f_{m})$ and $I_{n}(g_{n})$ are independent random variables for every $m,n\geq 1$. Using their joint characteristic function, it can easily be shown that $F$ and $G$ are independent (in the usual sense), as soon as they are strongly independent see also Proposition 8 in \cite{UZ}).

\begin{lemma}\label{ll3}
The two positive integers $q_1$ and $q_2$ are arbitrary.	Let $ F_1=I_{q_{1}}(f_{1})$ and $ F_{2}= I_{q_{2}} (f_{2}) $ with $f_{1} \in L ^{2} (\mathbb{R} ^{q_{1}})$ and $f_{2} \in L ^{2} (\mathbb{R} ^{q_{2}})$. Assume that $ F_{1}$ and $ F_{2}$ are independent. Then  for any integers $a,b \geq 1$, $ F_{1} ^{a}$ and $ F_{2}^{b}  $ are strongly independent.  
\end{lemma}
{\bf Proof: }  Without any restriction, one can suppose that the kernel functions $f_1$ and $f_2$, associated with $F_1$ and $F_2$, belong to $L^{2}_{S} ( \mathbb{R} ^{q_{1}})$ and $L^{2}_{S} ( \mathbb{R} ^{q_{2}})$ respectively; one mentions in passing that, for any integer $q\ge 0$, the closed subspace of $L^{2} ( \mathbb{R} ^{q})$ formed by the symmetric functions is denoted by $L^{2}_{S}( \mathbb{R} ^{q})$. 
Thus, in view of the criterion for the independence of multiple stochastic integrals in \cite{UZ}, the independence property of the random variables $F_1$ and $F_2$ can be reformulated as
\begin{equation}
	\label{22f-1}
	f_{1} \otimes _{1} f_{2} =0 \mbox { almost everywhere on } \mathbb{R} ^{q_{1} + q_{2}-2}.
\end{equation}

First we study the case $a=1$. We are going to prove, by induction on the integer $b\geq 1$, that the random variables $ F_{1}$  and $ F_{2} ^{b}$  are strongly independent.  For $b=1$, the claim is true by hypothesis. In the sequel, we assume, for any integer $b\ge 1$, that $ F_{1}$ and $ F_{2} ^{b}$ are strongly independent and we show that $ F_{1} $ and $ F_{2} ^{b+1}$ are strongly independent as well. We know from the product formula for multiple Wiener-It\^o integrals (\ref{prod}) that $ F_{2}^{b}$ admits the following chaos expansion:
\begin{equation}
\label{cha:f-2b}
	F_{2} ^{b}= \sum_{k=0} ^{bq_{2}} I_{k} (h_{2,k}) \mbox{ with } h_{2,k} \in L ^{2}_{S}(\mathbb{R} ^{k}). 
\end{equation}
Then,  the  induction hypothesis implies, for all $k=1,..., bq_{2}$, that
\begin{equation}
	\label{22f-2}
	f_{1} \otimes _{1} h_{2,k}=0 \mbox{ almost everywhere on } \mathbb{R} ^{q_{1}+ k-2}.
\end{equation}
Notice that, one can derive from \eqref{cha:f-2b} and the product formula for multiple Wiener-It\^o integrals that
\begin{eqnarray*}
	F_{2} ^{b+1} &=&\sum_{k=0} ^{bq_{2}} I_{k} (h_{2,k})I_{q_{2}}(f_{2})\\
	&=& \sum_{k=0} ^{bq_{2}} \sum_{r=0} ^{q_{2}\wedge k}r! \binom{k}{r} \binom{q_2}{r} I_{k+q_{2} -2r }\left( h_{2,k}\otimes _{r} f_{2}\right).
\end{eqnarray*}
It then suffices to prove that, for every $k=0,1,..., bq_{2}$ and $r=0,1,..., q_{2}\wedge k$ such that $k+q_{2}>2r$, the random  variables $F_{1}= I_{q_{1}} (f_{1}) $ and $ I_{k+q_{2} -2r }\left( h_{2,k}\otimes _{r} f_{2}\right)$  are independent. This is equivalent to 
\begin{equation}
	\label{22f-3}
	( h_{2,k} \tilde{\otimes }_r f_{2} ) \otimes _{1} f_{1}= 0 \mbox{ almost everywhere on } \mathbb{R} ^{q_{1}+ q_{2}+k -2r-2},
\end{equation}
where $h_{2,k} \tilde{\otimes }_r f_{2}$ denotes the symmetrization of $h_{2,k} \otimes_r f_{2}$.
We have,  by the definition of the contraction, 
\begin{eqnarray}
	&&	\left( ( h_{2,k} \tilde{\otimes }_r f_{2} ) \otimes _{1} f_{1}\right) (t_{1},...,t_{ k+q_{2}+q_{1}-2r-2})\nonumber \\
	&&= \int_{\mathbb{R}} dx ( h_{2,k} \tilde{\otimes }_r f_{2} )(t_{1},...,t_{k+q_{2} -2r-1}, x) f_{1} (t_{k+q_{2}-2r},..., t_{k+q_{2}-2r+q_{1}-2}, x).\label{22f-4}
\end{eqnarray}
We notice that for all $n\ge 2$ and for any function $ g\in L^{2}(\mathbb{R}^{n} )$, its symmetrization $\tilde{g}$, can be written, for every $(t_{1},...,t_{n-1},x)\in\mathbb{R}^{n}$, as 
\begin{equation*}
	\tilde{g} (t_{1},...,t_{n-1},x)= \frac{1}{ n!} \sum _{\sigma \in S_{n-1}}\sum_{i=1}^{n} g( t_{\sigma (1)},...,t_{\sigma (i-1)}, x, t_{\sigma (i)},...,t_{\sigma (n-1)}), 
\end{equation*}
where $S_{n-1}$ is the set of the permutations of $\{1,\ldots,n-1\}$. Thus,  by using the above formula in (\ref{22f-4}), we obtain that
\begin{eqnarray*}
	&&	\left( ( h_{2,k} \tilde{\otimes } f_{2} ) \otimes _{1} f_{1}\right) (t_{1},...,t_{ k+q_{2}+q_{1}-2r-2})\nonumber \\
	&&=\frac{1}{ (k+q_{2}-2r)!} \sum_{\sigma \in S_{k+q_{2}-2r-1}}\sum_{i=1} ^{k+q_{2}-2r} \int_{\mathbb{R}} dx \\
	&&\times (h_{2,k}\otimes _{r} f_{2}) (t_{\sigma(1)},...,t_{\sigma (i-1)}, x, t_{\sigma (i)},...,t_{\sigma (k+q_{2}-2r-1)}) \\
	&&\times f_{1} (t_{ k+q_{2}-2r},...,t_{k+q_{2}-2r+q_{1}-2}, x).
\end{eqnarray*}
Then, it turns out that in order to get the conclusion (\ref{22f-3}), it suffices to show that, for every $k=0,1,..., bq_{2}$ and $r=0,1,..., q_{2}\wedge k$ satisfying $k+q_{2}>2r$, and for all $\sigma \in S_{ k+q_{2}-2r-1}$ and $i=1,..., k+q_{2}-2r$, we have 
\begin{eqnarray}
&&	\int_{\mathbb{R}} dx (h_{2,k}\otimes _{r} f_{2}) (t_{\sigma(1)},...,t_{\sigma (i-1)}, x, t_{\sigma (i)},...,t_{\sigma (k+q_{2}-2r-1)})\nonumber \\
&&	\times f_{1} (t_{ k+q_{2}-2r},...,t_{k+q_{2}-2r+q_{1}-2}, x)=0\label{22f-5}
\end{eqnarray}
almost everywhere with respect to $t_{1},...,t_{k+q_{2}-2r+q_{1} -2}$.
Assume $k-r\leq i-1$. Then, using Fubini Theorem, the left-hand side of (\ref{22f-5}) writes 

	\begin{eqnarray*}
		&& \int_{\mathbb{R}} dx \int_{\mathbb{R}^r} d\,\overline{y}\times h_{2,k} ( t_{\sigma (1)},..., t_{\sigma (k-r)},\overline{y}) \\
		&&\times f_{2} (\overline{y},t_{\sigma (k-r+1)},...,t_{\sigma (i-1)}, x, t_{\sigma (i)},...,t_{\sigma (k+q_{2}-2r-1)})\\
		&&\times f_{1} (t_{ k+q_{2}-2r},...,t_{k+q_{2}-2r+q_{1}-2}, x)\\
		&&= \int_{\mathbb{R}^r} d\,\overline{y}\times h_{2,k} ( t_{\sigma (1)},..., t_{\sigma (k)}, \overline{y})\\
		&& \times\int_{\mathbb{R}} dx \times f_{2} (\overline{y},t_{\sigma (k+1)},...,t_{\sigma (i-1)}, x, t_{\sigma (i)},...,t_{\sigma (k+q_{2}-2r-1)})\\
		&&\times f_{1} (t_{ k+q_{2}-2r},...,t_{k+q_{2}-2r+q_{1}-2}, x),
	\end{eqnarray*}
	and the last integral on $\mathbb{R}$ vanishes due to (\ref{22f-1}).
Assume $k-r\geq i$. Then, using Fubini Theorem, the left-hand side of (\ref{22f-5}) can be written as 

	\begin{eqnarray*}
		&&\int_{\mathbb{R}^r} d\,\overline{y}\times f_{2} ( t_{\sigma(k-r)},...., t_{\sigma (k+q_{2}-2r-1)},\overline{y})\\
		&&\times \int_{\mathbb{R}}dx \times h_{2,k}(\overline{y},t_{\sigma(1)},...,t_{\sigma (i-1)},x, t_{\sigma(i)},..., t_{\sigma(k-r-1)})  f_{1} (t_{ k+q_{2}-2r},...,t_{k+q_{2}-2r+q_{1}-2}, x)
	\end{eqnarray*} 
and again the integral over $\mathbb{R}$ vanishes because of (\ref{22f-2}).

Let us now study the general case where the integer $a\geq 1$ is arbitrary. We are going to prove, by induction on $a$, that $F_{1}^{a}$ and $ F_{2}^{b}$ are strongly independent. We already have shown that the claim is true when $a=1$. Assume  $F_{1}^{a}$ and $ F_{2}^{b}$  are strongly independent and let us show that $F_{1}^{a+1}$ and $ F_{2}^{b}$ satisfy this same strong independence property.

It follows from the product formula for multiple Wiener-It\^o integrals (\ref{prod})  that $ F_{1}^{a}$ admits the following chaos expansion:
	\begin{equation}
\label{cha:f-1a}
		F_{1} ^{a}= \sum_{k=0}^{a q_{1}} I_{k}(h_{1,k}) \mbox{ with } h_{1,k} \in L ^{2}_{S}(\mathbb{R} ^{k}). 
	\end{equation}
Then, in view of the induction hypothesis, we have, for all integers $k\ge 1$ and $l\ge 1$,
\begin{equation}\label{22f-6}
	h_{1,k}\otimes_1 h_{2,l} =0 \mbox{ almost everywhere on } \mathbb{R} ^{k+l-2}.
\end{equation}

We also have from the case $a=1$, for every integer $l\ge 1$,
\begin{equation}\label{22f-7}
	f_{1}\otimes_1 h_{2,l} =0 \mbox{ almost everywhere on } \mathbb{R} ^{q_{1}+l-2}.
\end{equation}
Moreover, using \eqref{cha:f-1a} and the product formula (\ref{prod}), we get that
\begin{equation*}
	F_{1} ^{a+1}=  \sum_{k=0}^{a q_{1}}\sum_{r=0} ^{q_{1}\wedge k} r! \binom{q_1}{r} \binom{k}{r} I_{k+q_{1}-2r}(h_{1,k}\otimes _{r} f_{1}).
\end{equation*}
It then suffices to prove that, for all $l=1,..., bq_{2}$, $k=0,1,..., bq_{1}$ and $r=0,1,...,q_{1}\wedge k$ such that $k+q_{1}-2r>0$, the two random variables $ I_{l} ( h_{2,l})$ and $I_{k+q_{1}-2r}(h_{1,k}\otimes _{r} f_{1})$ are independent, which amounts to prove that

\begin{equation*}
	(h_{1,k}\tilde{\otimes} _{r} f_{1})\otimes _{1} h_{2,l}=0 \mbox{ almost everywhere on } \mathbb{R} ^{ k+q_{1}-2r+l-2}.
\end{equation*}

This follows exactly as in the case $a=1$, based on the hypotheses  (\ref{22f-6}) and (\ref{22f-7}). \qed

Before ending the present section, let us recall the following result from \cite{BDT} which is related to strongly independent random variables and Malliavin calculus (see Section \ref{app}).
\begin{lemma}\label{ll4}
	Let $ F,G$ be strongly independent random variables. Assume $F, G\in \mathbb{D}^{1,2}$. Then
	\begin{enumerate}
		\item  We have $ \langle DF, D(-L) ^{-1}G\rangle _{ L ^{2} (\mathbb{R})}= \langle DG, D(-L) ^{-1}F\rangle _{ L ^{2} (\mathbb{R})}=0.$
		\item The random variables $ \langle DF, D(-L) ^{-1}F \rangle _{ L ^{2} (\mathbb{R})}$ and $ \langle DG, D(-L) ^{-1}G \rangle _{ L ^{2} (\mathbb{R})}$ are strongly independent. 
	\end{enumerate}
\end{lemma}

\section{The modified power variation of the Hermite process}

Let $Z=(Z_{t}, t\geq 0) $ be a Hermite process and $p\geq 2$ an arbitrary  integer. The goal of this section is to analyse the asymptotic behavior of the modified power variation of $Z$ denoted by $S_{N, p}(Z)$. The difference between $S_{N, p}(Z)$ and usual power variation is that $S_{N, p}(Z)$ is defined by only considering the special increments of $Z$ introduced in Section \ref{sec22}. Namely, for each integer $N\ge 1$, it is defined as the following renormalized sum of the latter increments raised to the power $p$:

\begin{equation}
	\label{snpz}
	S_{N, p}(Z)=\frac{ 2 ^{pHN}}{ 2 ^{[N^{\gamma}]}}\sum _{l=1} ^{ 2 ^{[N ^{\gamma}]}}\left( \Delta Z_{l, N}\right) ^{p}.
\end{equation}

We first show that the modified power variation satisfies a law of large numbers.   
\begin{prop}
For all integer $p\geq 2$, one sets $\mu_{p}= \mathbf{E} Z_{1} ^{p}$. Let $S_{N, p}$ be given by (\ref{snpz}). Then 
	\begin{equation*}
		S_{N, p} (Z)\to _{N \to \infty} \mu_{p} \mbox{ in } L ^{1} (\Omega). 
	\end{equation*}
\end{prop}
{\bf Proof: } We use the decomposition (\ref{deco}), with $	\tilde{\Delta} Z_{l, N} $ and $ \check{\Delta} Z_{l, N} $ given by (\ref{dt}) and (\ref{dc}), respectively. In this way, we can write,
\begin{eqnarray}
	(\Delta Z_{l,N}) ^{p}&=& \sum_{j=0} ^{p} \binom{p}{j} (\tilde{\Delta}Z _{l,N})^{j} (\check{\Delta}Z_{l,N})^{p-j}\nonumber \\
	&=&(\tilde{\Delta} Z_{l,N})^{p} + \sum_{j=0} ^{p-1} \binom{p}{j} (\tilde{\Delta}Z _{l,N})^{j} (\check{\Delta}Z_{l,N})^{p-j},\label{8i-1}
\end{eqnarray}
and thus
\begin{eqnarray*}
	S_{N, p}(Z)- \mu_{p}&=& \frac{ 2 ^{pHN}}{ 2 ^{[N ^{\gamma}]}}\sum _{l=1} ^{ 2 ^{[N ^{\gamma}]}} \left( (\Delta Z_{l, N}) ^{p}- 2 ^{-pHN} \mu_{p} \right) \\
	&=& \frac{ 2 ^{pHN}}{ 2 ^{[N ^{\gamma}]}}\sum _{l=1} ^{ 2 ^{[N ^{\gamma}]}}\left( (\Delta Z_{l,N})^{p} - \mathbf{E}  (\Delta Z_{l,N})^{p} \right)\\
	&=& S_{N, p} ^{(1)} (Z) + (R_{N, p}(Z)-\mathbf{E}  R_{N, p}(Z)),
\end{eqnarray*}
with
\begin{equation*}
	S_{N, p} ^{(1)} (Z)=  \frac{ 2 ^{pHN}}{ 2 ^{[N ^{\gamma}]}}\sum _{l=1} ^{ 2 ^{[N ^{\gamma}]}}\left( (\tilde{\Delta} Z_{l,N})^{p} - \mathbf{E}  (\tilde{\Delta} Z_{l,N})^{p} \right)
\end{equation*}
and
\begin{equation}\label{rnpz}
	R_{N, p} (Z)= \frac{ 2 ^{pHN}}{ 2 ^{[N ^{\gamma}]}}\sum _{l=1} ^{ 2 ^{[N ^{\gamma}]}} \sum _{j=0} ^{p-1} \binom{p}{j} (\tilde{\Delta}Z _{l, N} )^{j} (\check{\Delta}Z _{l, N})^{p-j}.
\end{equation}
We first deal with the rest term $ R_{N, p}(Z)$. It follows from Cauchy-Schwarz's inequality, the estimates in Lemmas \ref{ll1} and \ref{ll2}, and \eqref{hyper} that
\begin{eqnarray}
	\mathbf{E} \vert R_{N, p} (Z)\vert &\leq &  \frac{ 2 ^{pHN}}{ 2 ^{[N ^{\gamma}]}}\sum _{l=1} ^{ 2 ^{[N ^{\gamma}]}}\sum _{j=0} ^{p-1} \binom{p}{j} \left( \mathbf{E}  (\tilde{\Delta}Z _{l, N} )^{2j} \right) ^{\frac{1}{2}} \left( \mathbf{E}  (\check{\Delta}Z _{l, N})^{2(p-j)}\right) ^{\frac{1}{2}}\nonumber\\
	&\leq & C 2 ^{ [N ^{\beta}]\frac{H-1}{q}}.\label{8i-2}
\end{eqnarray}
Consequently, $R_{N, p}(Z)$ converges to zero in $ L^{1}(\Omega)$ as $ N \to \infty$.  Next, we deal with the summand $ S_{N, p} ^{(1)} (Z)$ and we show that it converges to zero in $ L^{2}(\Omega)$ as $ N \to \infty$.  We have
\begin{eqnarray*}
	\mathbf{E}  \left( S_{N, p}^{(1) } (Z) \right) ^{2} &=&  \frac{ 2 ^{2pHN}}{ 2 ^{2[N ^{\gamma}]}}\sum _{l, l'=1} ^{ 2 ^{[N ^{\gamma}]}}\mathbf{E}  \left( (\tilde{\Delta} Z_{l,N})^{p} - \mathbf{E}  (\tilde{\Delta} Z_{l,N})^{p} \right)\left( (\tilde{\Delta} Z_{l',N})^{p} - \mathbf{E}  (\tilde{\Delta} Z_{l',N})^{p} \right)\\
	&=&  \frac{ 2 ^{2pHN}}{ 2 ^{2[N ^{\gamma}]}}\sum _{l=1} ^{ 2 ^{[N ^{\gamma}]}}\mathbf{E}  \left( (\tilde{\Delta} Z_{l,N})^{p} - \mathbf{E}  (\tilde{\Delta} Z_{l,N})^{p} \right)^{2}\\
	&\leq & C \frac{ 2 ^{2pHN}}{ 2 ^{2[N ^{\gamma}]}}\sum _{l=1} ^{ 2 ^{[N ^{\gamma}]}}\mathbf{E}  (\tilde{\Delta} Z_{l,N})^{2p} \\
	&\leq &C \frac{ 1}{ 2 ^{[N ^{\gamma}]}},
\end{eqnarray*}
where we used the independence property from Proposition \ref{pp1}, point 1.  and the estimate (\ref{15i-1}) in Lemma \ref{ll2}. We obtain
\begin{equation*}
	\mathbf{E}  \vert S_{N, p} (Z) -\mu_{p}\vert \leq C \left(  2 ^{ [N ^{\beta}]\frac{H-1}{q}}+ 2 ^{-\frac{ [ N ^{\gamma}]}{2}}\right) \leq C 2 ^{-\frac{ [ N ^{\gamma}]}{2}}.
\end{equation*}
\qed

The next step is to show that, after a proper renormalization, the power variation sequence (\ref{snpz}) satisfies a Central Limit Theorem.  For any integer number $p\geq 2$, let us introduce the sequence $ (V_{N, p}(Z), N\geq 1) $ given by 
\begin{eqnarray}\label{vnp}
	V_{N, p}(Z) &=& \sqrt{ 2 ^{ [N ^{\gamma}]}}\left( S_{N, p}(Z)-\mu _{p}\right) \\
	&=&\frac{ 2 ^{pHN}}{\sqrt{2 ^{ [N ^{\gamma}]}}}\sum_{l=1} ^{ 2 ^{[N ^{\gamma}]}} \left[ (\Delta Z _{l,N}) ^{p} - \mathbf{E}  (\Delta Z _{l,N}) ^{p}\right]. \nonumber
\end{eqnarray}
By (\ref{8i-1}), \begin{eqnarray}
\label{deco:vnp}
	V_{N,p}(Z)&=&  \frac{ 2 ^{pHN}}{\sqrt{2 ^{ [N ^{\gamma}]}}}\sum_{l=1} ^{ 2 ^{[N ^{\gamma}]}} \left[ (\tilde{\Delta} Z _{l,N}) ^{p} - \mathbf{E}  (\tilde{\Delta} Z _{l,N}) ^{p}\right]\nonumber\\
	&&+  \frac{ 2 ^{pHN}}{\sqrt{2 ^{ [N ^{\gamma}]}}}\sum_{l=1} ^{ 2 ^{[N ^{\gamma}]}} \sum_{j=0} ^{p-1} \binom{p}{j} \left[ (\tilde{\Delta}Z _{l,N})^{j} (\check{\Delta}Z_{l,N})^{p-j}-\mathbf{E} (\tilde{\Delta}Z _{l,N})^{j} (\check{\Delta}Z_{l,N})^{p-j}\right]\nonumber\\
	&=& V_{N,p}^{(1)} + \sqrt{ 2 ^{[N ^{\gamma}]}}( R_{N,p}(Z)- \mathbf{E}  R_{N, p}(Z)),
\end{eqnarray}
with
\begin{equation}\label{vnp1}
	V_{N,p}^{(1)}=  \frac{ 2 ^{pHN}}{\sqrt{2 ^{ [N ^{\gamma}]}}}\sum_{l=1} ^{ 2 ^{[N ^{\gamma}]}} \left[ (\tilde{\Delta} Z _{l,N}) ^{p} - \mathbf{E}  (\tilde{\Delta} Z _{l,N}) ^{p}\right]
\end{equation}
and $R_{N, p}(Z)$ given by (\ref{rnpz}). By (\ref{8i-2}),  
\begin{equation}
	\label{15i-3}
\sqrt{ 2 ^{[N ^{\gamma}]}}\mathbf{E}  \vert  R_{N,p}(Z)- \mathbf{E}  R_{N, p}(Z)\vert \leq C  2 ^{N^{\beta} \frac{H-1}{q}+ \frac{ N ^{\gamma}}{2}}.
\end{equation}

We will now focus on the asymptotic behaviour of the sequence $\left( V_{N, p} ^{(1)}, N\geq 1\right)$. We will prove that this sequence satisfies a Central Limit Theorem. We start by the study of the asymptotic behavior of the  second moment. 

\begin{lemma} 
	\label{lem:vaV1}

	Let $p\geq 2$ be an integer number. Consider the sequence $\left( V_{N, p} ^{(1)}, N\geq 1\right)$ given by (\ref{vnp1}). Then 
	\begin{equation}\label{5i-1}
		\left| \mathbf{E}  (V_{N,p} ^{(1)})^{2} - m_{p} \right|  \leq C 2 ^{ [N ^{\beta}]\frac{H-1}{q}},
	\end{equation}
	where 
	\begin{equation}
		\label{mp}
		m_{p}= {\rm Var}(Z_1^p) =\mathbf{E}  Z_{1} ^{2p}- (\mathbf{E}  Z_{1} ^{p} )^{2}=\mu _{2p}- \mu _{p} ^{2}. 
	\end{equation}
	In particular,
	\begin{equation*}
		\mathbf{E}  (V_{N,p} ^{(1)})^{2} \to_{ N \to \infty} m_{p}.
	\end{equation*}
\end{lemma}
{\bf Proof: } For $N\geq 1$ and $p\geq 1$, we write 
\begin{eqnarray*}
	\mathbf{E} (V_{N, p} ^{(1)}) ^{2} &=& \frac{ 2 ^{2pHN}}{2 ^{[N ^{\gamma}]}}\sum_{l,k=1} ^{2 ^{[N ^{\gamma}]}}\mathbf{E}  \left[ (\tilde{\Delta} Z _{l,N}) ^{p} - \mathbf{E}  (\tilde{\Delta} Z _{l,N}) ^{p}\right]\left[ (\tilde{\Delta} Z _{k,N}) ^{p} - \mathbf{E}  (\tilde{\Delta} Z _{k,N}) ^{p}\right]\\
	&=& \frac{ 2 ^{2pHN}}{2 ^{[N ^{\gamma}]}}\sum_{l=1} ^{2 ^{[N ^{\gamma}]}}\mathbf{E}  \left[ (\tilde{\Delta} Z _{l,N}) ^{p} - \mathbf{E}  (\tilde{\Delta} Z _{l,N}) ^{p}\right]^{2},
\end{eqnarray*}	
where we used the fact that $\tilde{\Delta} Z _{l,N}, l=1,..., 2 ^{[N ^{\gamma}]}$ are independent  (Proposition \ref{pp1}, point 1.). By using now the fact that  the random variables $\tilde{\Delta} Z _{l,N}, l=1,..., 2 ^{[N ^{\gamma}]}$ are identically distributed, we obtain 
\begin{equation}
\label{e:add1}
	\mathbf{E} (V_{N, p} ^{(1)}) ^{2} =  2 ^{2pHN}\mathbf{E}  \left[ (\tilde{\Delta} Z _{1,N}) ^{p} - \mathbf{E}  (\tilde{\Delta} Z _{1,N}) ^{p}\right]^{2}.
\end{equation}
We next write 
\begin{eqnarray}
	\label{ad-a1}
	\mathbf{E}  ( V_{N,p}^{(1)} ) ^{2} -m_{p} &=& 2 ^{2pHN}\mathbf{E}  (\tilde{\Delta} Z_{1, N })^{2p} -\mathbf{E}  Z_{1} ^{2p} \\
	&&-\left[ \left( 2 ^{pHN}\mathbf{E}  (\tilde{\Delta} Z_{1, N })^{p}\right) ^{2}-  (\mathbf{E}  Z_{1} ^{p}) ^{2}\right].\nonumber
\end{eqnarray}
We can derive from Lemma \ref{ll2} that 
\begin{equation}
	\label{ad-a2}
	\left|  2 ^{2pHN}\mathbf{E}  (\tilde{\Delta} Z_{1, N })^{2p} -\mathbf{E}  Z_{1} ^{2p}\right| \leq C_{2p}\ 2 ^{ [N ^{\beta}] \frac{H-1}{q}}.
\end{equation}
Moreover, using again Lemma \ref{ll2} and the inequality 
\[
\vert  2 ^{pHN}  \mathbf{E}  (\tilde{\Delta} Z_{1, N })^{p}+ \mathbf{E}  Z_{1} ^{p}\vert \leq \vert 2 ^{pHN}  \mathbf{E}  (\tilde{\Delta} Z_{1, N })^{p}- \mathbf{E}  Z_{1} ^{p}\vert +2 \mathbf{E} \vert  Z_{1} \vert^{p},
\]
we obtain that
\begin{eqnarray}
	\label{ad-a3}
	&&\left |\left( 2 ^{pHN}  \mathbf{E}  (\tilde{\Delta} Z_{1, N })^{p}\right) ^{2}-  (\mathbf{E}  Z_{1} ^{p}) ^{2}\right|\nonumber\\
	&&\leq \vert 2 ^{pHN}  \mathbf{E}  (\tilde{\Delta} Z_{1, N })^{p}- \mathbf{E}  Z_{1} ^{p}\vert \times \vert  2 ^{pHN}  \mathbf{E}  (\tilde{\Delta} Z_{1, N })^{p}+ \mathbf{E}  Z_{1} ^{p}\vert\nonumber\\
	&&\leq C_{p}' \ 2 ^{[N ^{\beta}] \frac{H-1}{q}}.
\end{eqnarray}
Finally, putting together (\ref{ad-a1}), (\ref{ad-a2}) and (\ref{ad-a3}), it follows that (\ref{5i-1}) holds. \qed

Now, we need to define some distances between the probability distributions of random variables. We refer to \cite{NP-book}, Appendix C, for more details. Usually, the  distance between the laws of two real-valued random variables $F$ and $G$ is defined as
\begin{equation}
	\label{dw}
	d (F, G)= \sup_{h\in \mathcal{A}}\left| \mathbf{E}h(F)-\mathbf{E}h(G)\right|,
\end{equation}
where $\mathcal{A}$ is a class of functions satisfying $h(F), h(G)\in L ^{1} (\Omega)$ for every $h\in \mathcal{A}$. When $\mathcal{A}$ is the set of Lipschitz continuous functions $h:\mathbb{R} \to \mathbb{R}$ such that $\Vert h\Vert _{Lip} \leq 1$, where
$$\Vert h\Vert _{Lip}= \sup_{x, y\in \mathbb{R} , x\not=y} \frac{ \vert h(x)-h(y)\vert} {\vert x-y\vert },$$
then (\ref{dw}) gives the Wasserstein distance. When $\mathcal{A}$ is the set of indicator functions $\{ \mathbbm{1}_{(-\infty, z]},\\ z\in \mathbb{R}\}$ then (\ref{dw}) gives the Kolmogorov distance,  for $\mathcal{A}= \{ \mathbbm{1}_{B}, B \in \mathcal{B}(\mathbb{R})\}$ we have the total variation distance, while the choice of $\mathcal{A}$ to be the class of functions $h$ with $\Vert h\Vert _{Lip}+ \Vert h\Vert _{\infty}<\infty$  leads to the Fortet-Mourier distance.

Let us recall a classical result from Stein-Malliavin calculus (see Theorem 5.1.3 and Remark 5.1.4 in \cite{NP-book}). Below, $d$ could be any of the above distances (Kolmogorov, Total variation, Wasserstein or Fortet-Mourier). We refer to the Appendix for the definition of the Malliavin derivative $D$ and of the Ornstein-Uhlenbeck operator $L$ with respect to an isonormal process $(W(h), h\in \mathcal{H})$, where  ${\mathcal{H}}$ is a real separable infinite-dimensional Hilbert space. Notice that in our present article, we have $\mathcal{H}= L ^{2} (\mathbb{R})$.

\begin{theorem}
	\label{tt1}
	Let $F$ be a random variable belonging to a finite sum of Wiener chaoses such that $\mathbf{\mathbf{E} }F=0$ and $ \mathbf{\mathbf{E} } F ^{2}= \sigma ^{2}$. Let $\gamma>0$. Then
	\begin{equation*}
		d(F, N(0, \gamma ^{2}))\leq C \left( \sqrt{{\rm Var} \left( \langle DF, D(-L)^{-1}F\rangle _{\mathcal{H}} \right) }+ \vert \sigma ^{2}- \gamma ^{2} \vert \right). 
	\end{equation*}
\end{theorem}

By applying Theorem \ref{tt1}, we will show that the sequence given by (\ref{vnp1}) converges in distribution to a Gaussian law and we also estimate the Wasserstein distance corresponding to this limit theorem.

\begin{prop}\label{pp3}
	Consider the sequence $(V_{N,p}^{(1)}, N\geq 1)$ given by (\ref{vnp1}). Then 
	\begin{equation*}
		V_{N, p}^{(1) }\to _{N \to\infty} ^{(d)}N(0, m_{p}) 
	\end{equation*}
	and for $N$ large enough, 
	\begin{equation*}
		d\left( V_{N,p}^{(1)}, N (0, m_{p})\right)\leq C 2 ^{ -\frac{ [N ^{\gamma}]}{2}}.
	\end{equation*}
\end{prop}
{\bf Proof: }Throughout the proof we assume that $N\ge 2^{(1-\beta)^{-1}}$; there is no restriction to make this assumption.  Let us evaluate the quantity ${\rm Var} \left( \langle DV_{N,p}^{(1)}, D(-L) ^{-1}V_{N, p} ^{(1)}\rangle _{L ^{2}(\mathbb{R})}\right)$. We have
\begin{equation*}
	DV_{N,p}^{(1)} = \frac{ 2 ^{pHN}}{\sqrt{2 ^{ [N ^{\gamma}]}}}\sum_{l=1} ^{2 ^{[N ^{\gamma}]}}D  \left[ (\tilde{\Delta} Z _{l,N}) ^{p} - \mathbf{E}  (\tilde{\Delta} Z _{l,N}) ^{p}\right]\
\end{equation*}
and
\begin{equation*}
	D(-L) ^{-1} V_{N,p}^{(1)} = \frac{ 2 ^{pHN}}{\sqrt{2 ^{ [N ^{\gamma}]}}}\sum_{l=1} ^{2 ^{[N ^{\gamma}]}}D (-L) ^{-1} \left[ (\tilde{\Delta} Z _{l,N}) ^{p} - \mathbf{E}  (\tilde{\Delta} Z _{l,N}) ^{p}\right].
\end{equation*}
 We know from Lemma \ref{ll3} that, for all $k,l=1,..., 2 ^{[N ^{\gamma}]}$ with $k\not=l$, the random variables $ (\tilde{\Delta} Z _{l,N}) ^{p} $ and $ (\tilde{\Delta} Z _{k,N}) ^{p} $ are strongly independent.  Then 
\begin{eqnarray*}
	&&\langle DV_{N,p}^{(1)}, D(-L) ^{-1}V_{N, p} ^{(1)}\rangle _{L ^{2}(\mathbb{R})}\\
	&=& \frac{ 2 ^{2pHN}}{2 ^{ [N ^{\gamma}]}}\sum_{l,k=1} ^{ 2 ^{[N ^{\gamma}]}} \langle D  \left[ (\tilde{\Delta} Z _{l,N}) ^{p} - \mathbf{E}  (\tilde{\Delta} Z _{l,N}) ^{p}\right], D (-L) ^{-1} \left[ (\tilde{\Delta} Z _{k,N}) ^{p} - \mathbf{E}  (\tilde{\Delta} Z _{k,N}) ^{p}\right]\rangle  _{L ^{2}(\mathbb{R})}\\
	&=& \frac{ 2 ^{2pHN}}{2 ^{ [N ^{\gamma}]}}\sum_{l=1} ^{ 2 ^{[N ^{\gamma}]}} \langle D  \left[ (\tilde{\Delta} Z _{l,N}) ^{p} - \mathbf{E}  (\tilde{\Delta} Z _{l,N}) ^{p}\right], D (-L) ^{-1} \left[ (\tilde{\Delta} Z _{l,N}) ^{p} - \mathbf{E}  (\tilde{\Delta} Z _{l,N}) ^{p}\right]\rangle  _{L ^{2}(\mathbb{R})},
\end{eqnarray*}
where we applied Lemma \ref{ll4}, point 1. Thus
\begin{eqnarray*}
	&&\langle DV_{N,p}^{(1)}, D(-L) ^{-1}V_{N, p} ^{(1)}\rangle _{L ^{2}(\mathbb{R})}\\
	&=&\mathbf{E} \langle DV_{N,p}^{(1)}, D(-L) ^{-1}V_{N, p} ^{(1)}\rangle _{L ^{2}(\mathbb{R})}+  \frac{ 2 ^{2pHN}}{2 ^{ [N ^{\gamma}]}}\sum_{l=1} ^{ 2 ^{[N ^{\gamma}]}} ( F_{l}- \mathbf{E}  F_{l}),
\end{eqnarray*}
where we used the notation 
\begin{equation*}
	F_{l}=\langle D  \left[ (\tilde{\Delta} Z _{l,N}) ^{p} - \mathbf{E}  (\tilde{\Delta} Z _{l,N}) ^{p}\right], D (-L) ^{-1} \left[ (\tilde{\Delta} Z _{l,N}) ^{p} - \mathbf{E}  (\tilde{\Delta} Z _{l,N}) ^{p}\right]\rangle  _{L ^{2}(\mathbb{R})}.
\end{equation*}
By point 2. in Lemma \ref{ll4}, the random variables $ F_{l}, l=1,...,2 ^{[N ^{\gamma}]}$ are pairwise independent.  Consequently, 
\begin{eqnarray}
	&& {\rm Var} \left(\langle DV_{N,p}^{(1)}, D(-L) ^{-1}V_{N, p} ^{(1)}\rangle _{L ^{2}(\mathbb{R})}\right) \nonumber \\
	&=&  \mathbf{E}  \left(\langle DV_{N,p}^{(1)}, D(-L) ^{-1}V_{N, p} ^{(1)}\rangle _{L ^{2}(\mathbb{R})}-\mathbf{E} \langle DV_{N,p}^{(1)}, D(-L) ^{-1}V_{N, p} ^{(1)}\rangle _{L ^{2}(\mathbb{R})}\right) ^{2}\nonumber  \\
	&=& \frac{ 2 ^{4pHN}}{2 ^{ 2[N ^{\gamma}]}}\sum_{l,k=1} ^{ 2 ^{[N ^{\gamma}]}}\mathbf{E}  \left( (F_{l}-\mathbf{E} F_{l}) (F_{k}-\mathbf{E} F_{k})\right) \nonumber \\
	&=&\frac{ 2 ^{4pHN}}{2 ^{ 2[N ^{\gamma}]}}\sum_{l=1} ^{ 2 ^{[N ^{\gamma}]}}\mathbf{E}  \left( (F_{l}-\mathbf{E} F_{l}) \right) ^{2}. \label{20f-2}
\end{eqnarray}
Let us show that for every $ l=1,...,2 ^{[N ^{\gamma}]}$, 
\begin{equation}
	\label{20f-1}
	\mathbf{E}  F_{l} ^{2} \leq C 2 ^{-4pHN}.
\end{equation}
Using Cauchy-Schwarz's inequality twice,  we obtain that
\begin{eqnarray*}
	\mathbf{E}  F_{l} ^{2} &\leq & C \mathbf{E}  \langle DV_{N,p} ^{(1)} , D(-L) ^{-1} V_{N,p}^{(1)} \rangle ^{2} _{L ^{2}(\mathbb{R})}\\
	&\leq &  C \mathbf{E}  \Vert DV_{N,p}^{(1)}\Vert ^{2} _{L ^{2}(\mathbb{R})} \Vert D(-L) ^{-1} V_{N,p} ^{1} \Vert ^{2} _{L ^{2}(\mathbb{R})}\\
	&\leq & C \left( \mathbf{E}  \Vert DV_{N,p}^{(1)}\Vert ^{4} _{L ^{2}(\mathbb{R})}\right) ^{\frac{1}{2}} \left( \mathbf{E} \Vert D(-L) ^{-1} V_{N,p} ^{1} \Vert ^{4} _{L ^{2}(\mathbb{R})}\right) ^{\frac{1}{2}}.
\end{eqnarray*}
On the other hand, since $ V_{N, p} ^{(1)} $ is a random variable belonging to a finite sum of Wiener chaoses, 
\begin{equation*}
	\mathbf{E}   \Vert DV_{N,p}^{(1)}\Vert ^{2} _{L ^{2}(\mathbb{R})}  \leq  C \mathbf{E}  \vert V_{N,p} ^{(1)} \vert ^{2} \leq C 2 ^{-2pHN},
\end{equation*}

and
\begin{equation*}
\mathbf{E} 	\Vert D(-L) ^{-1} V_{N,p} ^{(1)} \Vert ^{2} _{L ^{2}(\mathbb{R})}\leq C \mathbf{E}  \vert V_{N,p} ^{(1)} \vert ^{2} \leq C 2 ^{-2pHN}.
\end{equation*}
Moreover, since $ \Vert DV_{N,p}^{(1)}\Vert ^{2} _{L ^{2}(\mathbb{R})} $ and $\Vert D(-L) ^{-1} V_{N,p} ^{(1)} \Vert ^{2} _{L ^{2}(\mathbb{R})}$ are also random variables in a finite sum of Wiener chaoses, we have via the hypercontractivity property (\ref{hyper})
\begin{equation*}
		\mathbf{E}   \Vert DV_{N,p}^{(1)}\Vert ^{4} _{L ^{2}(\mathbb{R})} \leq C \left( 	\mathbf{E}   \Vert DV_{N,p}^{(1)}\Vert ^{2} _{L ^{2}(\mathbb{R})} \right) ^{2}\leq C 2 ^{-4pHN}
\end{equation*}
and 
\begin{equation*}
	\mathbf{E} 	\Vert D(-L) ^{-1} V_{N,p} ^{(1)} \Vert ^{4} _{L ^{2}(\mathbb{R})}\leq C  \left( \mathbf{E} 	\Vert D(-L) ^{-1} V_{N,p} ^{(1)} \Vert ^{2} _{L ^{2}(\mathbb{R})}\right) ^{2}\leq C 2 ^{-4pHN}.
\end{equation*}
So, the bound   (\ref{20f-1}) is proven.  By plugging this inequality into (\ref{20f-2}), we get 
\begin{equation}\label{20f-3}
	{\rm Var} \left(\langle DV_{N,p}^{(1)}, D(-L) ^{-1}V_{N, p} ^{(1)}\rangle _{L ^{2}(\mathbb{R})}\right) \leq C 2 ^{-[N^{\gamma}]}.
\end{equation}
Finally, putting together (\ref{5i-1}), (\ref{20f-3}) and Theorem \ref{tt1} we obtain the proposition. \qed

We can now conclude the convergence of the modified $p$-variation of the Hermite process. 

\begin{theorem}
	Consider the sequence $(V_{N,p}(Z), N\geq 1)$ given by (\ref{vnp}). Then 
	\begin{equation*}
		V_{N, p} (Z)\to ^{(d)} _{N \to \infty} N (0, m_{p}),
	\end{equation*}
	where $m_{p}={\rm Var}(Z_1^p) $ (see (\ref{mp})).  Moreover, denoting by $d_W$ the Wasserstein distance, for $N$ sufficiently large, we have
	\begin{equation*}
		d_{W} \left( V_{N,p}(Z), N(0, m_{p})\right)  \leq C 2 ^{-\frac{N ^{\gamma}}{2}}
	\end{equation*}
\end{theorem}
{\bf Proof: }The convergence in distribution follows from \eqref{deco:vnp}, Proposition \ref{pp3} and \eqref{15i-3} which implies that the random variable $\sqrt{ 2 ^{[N ^{\gamma}]}}( R_{N,p}(Z)- \mathbf{E}  R_{N, p}(Z))$ converges to $0$ in $L^1(\Omega)$ when $N \to \infty$. To get the estimate of the Wasserstein distance, we first notice that 
$$d_{W}( V_{n,p}(Z), V_{N,p} ^{(1)})\leq \mathbf{E}  \vert V_{N,p}(Z) -V_{N,p}^{(1)} \vert = \mathbf{E}  \sqrt{ 2 ^{[N ^{\gamma}]}}\vert  R_{N,p}(Z)- \mathbf{E}  R_{N, p}(Z)\vert ,$$
with $ R_{N,p}(Z)$ given by (\ref{rnpz}). Thus, for $N$ large enough,
\begin{eqnarray*}
	d_{W} ( V_{N,p}, N(0,m_{p}))&\leq & d_{W} ( V_{N,p}^{(1)}, N(0,m_{p})) + d_{W} ( V_{n,p}(Z), V_{N,p} ^{(1)})\\&\leq &  d_{W} ( V_{N,p}^{(1)}, N(0,m_{p})) + \sqrt{ 2 ^{[N ^{\gamma}]}}\mathbf{E}  \vert R_{N,p}\vert \leq  C 2 ^{-\frac{N ^{\gamma}}{2}}, 
\end{eqnarray*}
due to Proposition \ref{pp3} and the bound (\ref{8i-2}). 
\qed

\section{Weighted power variation and application to volatility estimation}

In this part, the purpose is to extend the results from the previous section, by looking to the asymptotic behavior of the weighted (modified) power variation of the Hermite process. In the second part of this section, we also apply the results to the volatility estimation in a simple model with Hermite noise. 

\subsection{Weighted power variation of the Hermite process}

We start by introducing the weighted modified power variation of the Hermite process. Let $ h:[0,1]\to \mathbb{R}$  be a measurable deterministic function. For every integer $p\geq 2$, we consider the sequence $(U_{N, p}(Z), N\geq 1)$  defined by
\begin{equation}\label{unpz}
	U_{N, p}(Z) = \frac{ 2 ^{pHN}}{\sqrt{2 ^{[N ^{\gamma}]}}} \sum _{l=1} ^{2 ^{ [N ^{\gamma}]}}h\left( \frac{l}{ 2 ^{[N ^{\gamma}]}}\right) \left( (\Delta Z_{l, N})^{p}  -\mathbf{E} (\Delta Z_{l, N})^{p}\right),
\end{equation}
where, for each $l=1,..., 2 ^{[{N ^{\gamma}]}}$, $\Delta Z_{l, N}$ is the special increment of the Hermite process given by (\ref{dz}). We study the limit behavior in distribution, as $N \to \infty$, of the sequence (\ref{unpz}).

Let  $(W_{t}, t\in [0,1])$ be a Brownian motion independent from the Wiener process $B$ in (\ref{hermite1}). For any deterministic function $g:[0,1]\to \mathbb{R}$, $g\in L^{2}([0,1])$, we let
\begin{equation*}
	\mathcal{I}(g)=\int_{0}^{1} g(s) dW_{s},
\end{equation*}
be the Wiener integral of $g$ with respect to $W$.

\begin{theorem}\label{tt3bi}
	Let $p\geq 2$ be an integer number and let $m_{p}$ be as in (\ref{mp}).  We assume that the function $h$ in (\ref{unpz}) satisfies the following assumption:
	\begin{equation}
		\label{h}
h \mbox{ is (bounded) and $\alpha$-H\"older continuous on the interval $[0,1]$ with $\alpha\in (0,1)$.}
\end{equation}

 Then, for each $\nu\in (0,\gamma)$, there is a constant $C$, which depends on $\nu$ and $h$, such that, for all $N$ large enough, one has 
	\begin{equation}
	\label{tt3bi:e1}
		d_W(U_{N, p} (Z), \sqrt{m_{p}}\,\mathcal{I}(h))\le C 2 ^{-\alpha N ^{\nu}}.
	\end{equation}
In particular, 
\begin{equation*}
	U_{N, p} (Z)\to ^{(d)} _{N \to \infty} \sqrt{m_{p}}\,\mathcal{I}(h).
\end{equation*}
	 	\end{theorem}

The proof of Theorem \ref{tt3bi} mainly relies on the following general lemma on Wasserstein distance.

\begin{lemma}
\label{lem:was}
Let $h:\mathbb{R}\rightarrow\mathbb{R}$ be an arbitrary measurable bounded deterministic function and $\|h\|_\infty=\sup_{v\in\mathbb{R}}|h(v)|$. Let $U$, $V$, $X$ and $Y$ be four real-valued random variables of $L^2(\Omega)$ such that the $2$-dimensional random vectors $(U,V)$ and $(X,Y)$ are independent. Then, one has
\begin{equation}
\label{lem:was:e1}
d_W\big (U+h(V)X,U+h(V)Y\big)\le \|h\|_\infty d_W (X,Y).
\end{equation}
\end{lemma}

\noindent {\bf Proof: } Denote by $\mathcal{A}_W$ the set of the Lipschitz continuous functions $f:\mathbb{R} \to \mathbb{R}$ such that $\Vert f\Vert _{Lip} \leq 1$. Observe that for any such $f$, one has, for all $z\in \mathbb{R}$, $|f(z)|\le |z|+|f(0)|$. the latter fact combined with the assumptions $U, V, X, Y\in L^2(\Omega)$ and $\|h\|_\infty<\infty$ imply that the expectation $\mathbf{E}\big (
f(U+h(V)X)-f(U+h(V)Y)\big)$ is well-defined and finite. Next one denotes $\mathbf{P}_{U,V}$ and  $\mathbf{P}_{X,Y}$ the probability distributions of the two random vectors $(U,V)$ and $(X,Y)$. Using the fact that these two random vectors are independent and Fubini Theorem one gets that 
\begin{eqnarray}
\label{lem:was:e2}
&&\Big |\mathbf{E}\big (
f(U+h(V)X)-f(U+h(V)Y)\big)\Big | \nonumber\\
&&=\Big | \int_{\mathbb{R}^2} d \mathbf{P}_{U,V}(u,v)\int_{\mathbb{R}^2} \Big (f(u+h(v)x)-f(u+h(v)y)\Big) d \mathbf{P}_{X,Y}(x,y)\Big|\nonumber\\
&&\le \int_{\mathbb{R}^2} d \mathbf{P}_{U,V}(u,v)\Big|\int_{\mathbb{R}^2} \Big (f(u+h(v)x)-f(u+h(v)y)\Big) d \mathbf{P}_{X,Y}(x,y)\Big| \nonumber\\
&& =\int_{\mathbb{R}^2} d \mathbf{P}_{U,V}(u,v)\Big | \mathbf{E}\Big (f(u+h(v)X)-f(u+h(v)Y)\Big)\Big|\nonumber\\
&&=\|h\|_\infty \int_{\mathbb{R}^2} d \mathbf{P}_{U,V}(u,v) \big | \mathbf{E}\big (g_{u,v}(X)-g_{u,v}(Y)\big)\big|,
\end{eqnarray}
where, for each fixed $(u,v)\in\mathbb{R}^2$, $g_{u,v}$ denotes the function of $\mathcal{A}_W$ defined, for every $z\in\R$, as 
$g_{u,v}(z)=\|h\|_\infty^{-1}f(u+h(v)z)$. Then, one can derive \eqref{lem:was:e2} and \eqref{dw} with $\mathcal{A}=\mathcal{A}_W$ that
\begin{eqnarray*}
&&\Big |\mathbf{E}\big (f(U+h(V)X)-f(U+h(V)Y)\big)\Big | \\
&&\le \|h\|_\infty \int_{\mathbb{R}^2} d \mathbf{P}_{U,V}(u,v) d_W (X,Y)= \|h\|_\infty  d_W (X,Y),
\end{eqnarray*}
and consequently that 
\begin{eqnarray*}
d_W\big (U+h(V)X, U+h(V)Y\big)&=&\sup_{f\in \mathcal{A}_W} \Big |\mathbf{E}\big (
f(U+h(V)X)-f(U+h(V)Y)\big)\Big |\\
&& \le \|h\|_\infty  d_W (X,Y).
\end{eqnarray*}
\qed
\\
\\ 
{\bf Proof of Theorem \ref{tt3bi}: } Throughout the proof, we fix the arbitrary integer  $N\ge 2^{(1-\beta)^{-1}}$. Let $U_{N}=U_{N, p} (Z)$ and $h_N:[0,1)\rightarrow \mathbb{R}$ be the step function which is, for all $k\in\{1,\ldots, 2 ^{[N^{\nu}]}\}$, equal to $h\left( \frac{k-1}{2 ^{[N^{\nu}]}}\right)$ on the interval $\left [\frac{k-1}{2 ^{[N^{\nu}]}}, \frac{k}{2 ^{[N^{\nu}]}}\right)$. First, notice that using the triangle inequality, one gets that 
\begin{eqnarray}
		\label{tt3bi:e30}
&& d_W(U_{N}, \sqrt{m_{p}}\,\mathcal{I}(h))\nonumber\\
&&\le d_W(U_{N}, \widetilde{U}_{N})+d_W(\widetilde{U}_{N}, \widetilde{B}_{N})+d_W(\widetilde{B}_{N}, \sqrt{m_{p}}\,\mathcal{I}(h_N))+d_W(\sqrt{m_{p}}\,\mathcal{I}(h_N), \sqrt{m_{p}}\,\mathcal{I}(h)). \nonumber\\
&& \le \mathbf{E}|U_{N}-\widetilde{U}_{N}|+\mathbf{E}|\widetilde{U}_{N}-\widetilde{B}_{N}|+d_W(\widetilde{B}_{N}, \sqrt{m_{p}}\,\mathcal{I}(h_N))+\sqrt{m_{p}}\, \mathbf{E}| \mathcal{I}(h_N)-\mathcal{I}(h)|.
\end{eqnarray}		
In the above inequality we have used the following notations: 
\begin{equation}
\label{tt3bi:e31}
\widetilde{U}_{N}=\frac{ 2 ^{pHN}}{\sqrt{2 ^{[N ^{\gamma}]}}} \sum _{l=1} ^{2 ^{ [N ^{\gamma}]}}h\left( \frac{l}{ 2 ^{[N ^{\gamma}]}}\right) \left( (\tilde{\Delta} Z_{l, N})^{p}  -\mathbf{E} (\tilde{\Delta} Z_{l, N})^{p}\right);
\end{equation}
\begin{equation}\label{tt3bi:e4b}
\widetilde{B}_N=\sum_{k=1}^{2^{[N^{\nu}]}} h\left( \frac{k-1}{2 ^{[N^{\nu}]}}\right) \widetilde{v}_{N} ^{(k)},
\end{equation}	
where, for every $k\in\{1, \ldots, 2^{[N^{\nu}]}\}$, 
\begin{equation}\label{tt3bi:e3b}
		\widetilde{v}_{N} ^{(k)} =\frac{ 2 ^{pHN}}{\sqrt{2 ^{ [N ^{\gamma}]}}}\sum_{l=(k-1) \delta(N)+1}^{l=k\delta (N)} \left[ (\tilde{\Delta} Z _{l,N}) ^{p} - \mathbf{E}  (\tilde{\Delta} Z _{l,N}) ^{p}\right],
	\end{equation}	
the positive integer $\delta (N)$ being defined as 
		\begin{equation}
		\label{tt3bi:e2}
		\delta (N)= 2^{[N^\gamma]-[N^{\nu}]};
		\end{equation}
and
\begin{equation}\label{tt3bi:e4t}
\mathcal{I}(h_N)=\sum_{k=1}^{2^{[N^{\nu}]}} h\left( \frac{k-1}{2 ^{[N^{\nu}]}}\right) \Delta^{(\nu)} W _{k-1,N},
\end{equation}	
where, for all $k\in\{1, \ldots, 2^{[N^{\nu}]}\}$,
\begin{equation}\label{tt3bi:e5}
 \Delta^{(\nu)} W _{k-1,N}=W_{\frac{k}{2^{[N^{\nu}]}}}-W_{\frac{k-1}{2^{[N^{\nu}]}}}.
\end{equation}

From now on, our goal is to provide an appropriate upper bound for each one of the four terms in the right-hand side of \eqref{tt3bi:e30}.

Since $h$ is a bounded function, combining \eqref{unpz} and \eqref{tt3bi:e31}, one obtain that 
\[
\mathbf{E}|U_{N}-\widetilde{U}_{N}|\le \frac{ 2 ^{pHN+1}\|h\|_\infty}{\sqrt{2 ^{[N ^{\gamma}]}}} \sum _{l=1} ^{2 ^{ [N ^{\gamma}]}}\mathbf{E}\left | (\Delta Z_{l, N})^{p}-(\tilde{\Delta} Z_{l, N})^{p}\right|.
\]
Then, one can derive from \eqref{8i-1}, \eqref{rnpz} and \eqref{8i-2} that
\begin{equation}
\label{tt3bi:e32}
\mathbf{E}|U_{N}-\widetilde{U}_{N}|\le C  2 ^{N^{\beta} \frac{H-1}{q}+ \frac{ N ^{\gamma}}{2}}.
\end{equation}

Let us now bound $\mathbf{E}|\widetilde{U}_{N}-\widetilde{B}_{N}|$. Observe that, one can derive from \eqref{tt3bi:e31}, \eqref{tt3bi:e4b} and \eqref{tt3bi:e3b} that
\begin{equation}
\label{tt3bi:e33}
\widetilde{U}_{N}-\widetilde{B}_{N}=\frac{ 2 ^{pHN}}{\sqrt{2 ^{ [N ^{\gamma}]}}}\sum_{k=1}^{2^{[N^{\nu}]}}\sum_{l=(k-1) \delta(N)+1}^{l=k\delta (N)} \left( h\left( \frac{l}{2 ^{[N^{\gamma}]}}\right) - h\left( \frac{k-1}{2 ^{[N^{\nu}]}}\right) \right)  \left[ (\tilde{\Delta} Z _{l,N}) ^{p} - \mathbf{E}  (\tilde{\Delta} Z _{l,N}) ^{p}\right].
\end{equation}
Also, observe that, for all $k\in\{1,\ldots, 2^{[N^\nu]}\}$ and $l\in\{(k-1)\delta(N)+1,\ldots k\delta(N)\}$, using \eqref{tt3bi:e2}, one has that
\[
\Big | \frac{l}{2 ^{[N^{\gamma}]}}-\frac{k-1}{2 ^{[N^{\nu}]}}\Big|=\frac{l-(k-1)\delta(N)}{2 ^{[N^{\gamma}]}}\le \frac{k\delta(N)-(k-1)\delta(N)}{2 ^{[N^{\gamma}]}}=2^{-[N^\nu]}.
\]
Then, the fact that $h$ is an $\alpha$-H\"older continuous function on the interval $[0,1]$ implies that 
\begin{equation}
\label{tt3bi:e34}
\left | h\left( \frac{l}{2 ^{[N^{\gamma}]}}\right) - h\left( \frac{k-1}{2 ^{[N^{\nu}]}}\right) \right |\le C 2^{-\alpha [N^\nu]}.
\end{equation} 
Next, one can derive from \eqref{tt3bi:e33},  point 1. in Proposition \ref{pp1}, \eqref{tt3bi:e34} and \eqref{15i-1} that
\begin{eqnarray*}
	\mathbf{E}|\widetilde{U}_{N}-\widetilde{B}_{N}|^{2} &=& \frac{ 2 ^{2pHN}}{2 ^{ [N ^{\gamma}]}}\sum_{k=1}^{2^{[N^{\nu}]}}\sum_{l=(k-1) \delta(N)+1}^{l=k\delta (N)} \left( h\left( \frac{l}{2 ^{[N^{\gamma}]}}\right) - h\left( \frac{k-1}{2 ^{[N^{\nu}]}}\right) \right)^{2}\\
	&&\times  \mathbf{E} \left[ \left( (\tilde{\Delta} Z _{l,N}) ^{p} - \mathbf{E}  (\tilde{\Delta} Z _{l,N}) ^{p}\right) ^{2} \right]\\
	&\leq & C 2^{-2\alpha [N^\nu]} \frac{ 2 ^{2pHN}}{2 ^{ [N ^{\gamma}]}}\sum_{k=1}^{2^{[N^{\nu}]}}\sum_{l=(k-1) \delta(N)+1}^{l=k\delta (N)}\mathbf{E} (\tilde{\Delta} Z _{l,N}) ^{2p}\\
	&\leq & C 2^{-2\alpha [N^\nu]}.
\end{eqnarray*}
Thus, using Cauchy-Schwarz's inequality one gets that
\begin{equation}
\label{tt3bi:e35}
\mathbf{E}|\widetilde{U}_{N}-\widetilde{B}_{N}|\le \sqrt{\mathbf{E}|\widetilde{U}_{N}-\widetilde{B}_{N}|^2}\le C 2^{-\alpha [N^\nu]}.
\end{equation}

Let us now bound $\mathbf{E}| \mathcal{I}(h_N)-\mathcal{I}(h)|$. It follows from Cauchy-Schwarz's inequality, the isometry property of the Wiener integral $\mathcal{I}$, the definition of the step function $h_N$ and the fact that $h$ is an $\alpha$-H\"older continuous function on the interval $[0,1]$ that 
\begin{eqnarray*}
&& \mathbf{E}| \mathcal{I}(h_N)-\mathcal{I}(h)|^2=\int_0^1 |h_N (s)-h(s)|^2 ds\\
&&=\sum_{k=1}^{2^{[N^{\nu}]}}\int_{\frac{k-1}{2 ^{[N^{\nu}]}}}^{\frac{k}{2 ^{[N^{\nu}]}}}\Big | h\Big( \frac{k-1}{2 ^{[N^{\nu}]}}\Big)-h(s)\Big |^2 ds\le C 2^{-2\alpha [N^\nu]}.
\end{eqnarray*}
Thus, using Cauchy-Schwarz's inequality one gets that
\begin{equation}
\label{tt3bi:e36}
\mathbf{E}| \mathcal{I}(h_N)-\mathcal{I}(h)|\le \sqrt{\mathbf{E}| \mathcal{I}(h_N)-\mathcal{I}(h)|^2}\le C 2^{-\alpha [N^\nu]}.
\end{equation}

Let us now bound $d_W (\widetilde{B}_N, \sqrt{m_{p}}\,\mathcal{I}(h_N))$. To this end, one needs to introduce some additional notations. For all $L\in \{0,\ldots, 2^{[N^{\nu}]}\}$, one sets 
\begin{equation}
\label{tt3bi:e7}
\widetilde{B}_{N,L}=\sum_{k=1}^{L} h\left( \frac{k-1}{2 ^{[N^{\nu}]}}\right) \widetilde{v}_{N} ^{(k)}+\sqrt{m_{p}}\sum_{k=L+1}^{2^{[N^{\nu}]}}h\left( \frac{k-1}{2 ^{[N^{\nu}]}}\right) \Delta^{(\nu)} W _{k-1,N}.
\end{equation}
Moreover, one sets
\begin{equation}
\label{tt3bi:e7a}
\widetilde{B}_{N,0}'=\sqrt{m_{p}}\sum_{k=2}^{2^{[N^{\nu}]}}h\left( \frac{k-1}{2 ^{[N^{\nu}]}}\right) \Delta^{(\nu)} W _{k-1,N},
\end{equation}
and, for every $L\in \{1,\ldots, 2^{[N^{\nu}]}\}$,
\begin{equation}
\label{tt3bi:e7b}
\widetilde{B}_{N,L}'=\sum_{k=1}^{L-1} h\left( \frac{k-1}{2 ^{[N^{\nu}]}}\right) \widetilde{v}_{N} ^{(k)}+\sqrt{m_{p}}\sum_{k=L+1}^{2^{[N^{\nu}]}}h\left( \frac{k-1}{2 ^{[N^{\nu}]}}\right) \Delta^{(\nu)} W _{k-1,N}.
\end{equation}
It easily follows from \eqref{tt3bi:e7}, \eqref{tt3bi:e7a} and \eqref{tt3bi:e7b} that, for all $L\in \{1,\ldots, 2^{[N^{\nu}]}\}$,
\begin{equation}
\label{tt3bi:e7c}
\widetilde{B}_{N,L}=\widetilde{B}_{N,L}'+h\left( \frac{L-1}{2 ^{[N^{\nu}]}}\right) \widetilde{v}_{N} ^{(L)}
\end{equation}
and 
\begin{equation}
\label{tt3bi:e7d}
\widetilde{B}_{N,L-1}=\widetilde{B}_{N,L}'+\sqrt{m_{p}}\, h\left( \frac{L-1}{2 ^{[N^{\nu}]}}\right) \Delta^{(\nu)} W _{L-1,N}.
\end{equation}
Moreover, the fact that the random variables $\widetilde{v}_{N} ^{(1)}, \ldots, \widetilde{v}_{N} ^{(2^{[N^{\nu}]})}, \Delta^{(\nu)} W _{0,N},\ldots \Delta^{(\nu)} W _{2^{[N^{\nu}]},N}$ are mutually independent implies that, for every $L\in \{1,\ldots, 2^{[N^{\nu}]}\}$, the 
random variable $\widetilde{B}_{N,L}'$ and the $2$-dimensional random vector $\big (\sqrt{m_{p}}\, \Delta^{(\nu)} W _{L-1,N}, \widetilde{v}_{N} ^{(L)}\big)$ are independent as well. Next, observe that, one can derive from \eqref{tt3bi:e4b}, \eqref{tt3bi:e4t}, \eqref{tt3bi:e7}, the triangle inequality, \eqref{tt3bi:e7c} and \eqref{tt3bi:e7d} that
\begin{eqnarray*}
&& d_W (\widetilde{B}_N, \sqrt{m_{p}}\,\mathcal{I}(h_N))=d_W (\widetilde{B}_{N,2^{[N^\nu]}}, \widetilde{B}_{N,0}) \le \sum_{L=1}^{2^{[N^\nu]}} d_W (\widetilde{B}_{N,L-1}, \widetilde{B}_{N,L})\\
&&=\sum_{L=1}^{2^{[N^\nu]}} d_W \left ( \widetilde{B}_{N,L}'+\sqrt{m_{p}}\, h\left( \frac{L-1}{2 ^{[N^{\nu}]}}\right) \Delta^{(\nu)} W _{L-1,N}, \widetilde{B}_{N,L}'+h\left( \frac{L-1}{2 ^{[N^{\nu}]}}\right) \widetilde{v}_{N} ^{(L)} \right).
\end{eqnarray*}
Thus, it results from Lemma \ref{lem:was} that 
\begin{equation}
\label{tt3bi:e8}
 d_W (\widetilde{B}_N, \sqrt{m_{p}}\,\mathcal{I}(h_N))\le \| h\|_\infty \sum_{L=1}^{2^{[N^\nu]}} d_W \left ( \sqrt{m_{p}}\, \Delta^{(\nu)} W _{L-1,N}\,, \widetilde{v}_{N} ^{(L)} \right).
 \end{equation}
 Next, observe that it follows from \eqref{tt3bi:e3b} and point 1. in Proposition \ref{pp1} that the random variables $\widetilde{v}_{N} ^{(L)}$, $L\in\{1,\ldots 2^{[N^\nu]}\}$, are identically distributed. Also observe that \eqref{tt3bi:e5} implies that the random variables $\sqrt{m_{p}}\, \Delta^{(\nu)} W _{L-1,N}$, $L\in\{1,\ldots 2^{[N^\nu]}\}$, are identically distributed. Thus, since, for every 
$L\in\{1,\ldots 2^{[N^\nu]}\}$, the two random variables $\widetilde{v}_{N} ^{(L)}$ and $\sqrt{m_{p}}\, \Delta^{(\nu)} W _{L-1,N}$ are independent, it turns out that the random vectors $\left ( \sqrt{m_{p}}\, \Delta^{(\nu)} W _{L-1,N}\,, \widetilde{v}_{N} ^{(L)} \right)$, $L\in\{1,\ldots 2^{[N^\nu]}\}$, are identically distributed. Then, one can derive from \eqref{tt3bi:e8} and the \\ equality $ \sqrt{m_{p}}\, \Delta^{(\nu)} W _{0,N} = \sqrt{m_{p}}\,W_{\frac{1}{2^{[N^\nu]}}}$ that
\begin{equation}
\label{tt3bi:e9}
 d_W (\widetilde{B}_N, \sqrt{m_{p}}\,\mathcal{I}(h_N))\le \| h\|_\infty 2^{[N^\nu]} d_W \left (  \sqrt{m_{p}}\,W_{\frac{1}{2^{[N^\nu]}}}\,, \widetilde{v}_{N} ^{(1)} \right).
\end{equation}
Since $\sqrt{m_{p}}\,W_{\frac{1}{2^{[N^\nu]}}}$ is a centred Gaussian random variable, in order to bound the Wassertein distance $d_W \left (  \sqrt{m_{p}}\,W_{\frac{1}{2^{[N^\nu]}}}\,, \widetilde{v}_{N} ^{(1)} \right)$ we will make use of the fundamental Theorem \ref{tt1}. Let us first bound $\Big |{\rm Var}\big (\sqrt{m_{p}}\,W_{\frac{1}{2^{[N^\nu]}}}\big)-{\rm Var}\big (\widetilde{v}_{N} ^{(1)}\big)\Big |$. One clearly has that 
\begin{equation}
\label{tt3bi:e10}
{\rm Var}\big (\sqrt{m_{p}}\,W_{\frac{1}{2^{[N^\nu]}}}\big)=m_p \, 2^{-[N^\nu]}.
\end{equation}
Moreover, in view of \eqref{tt3bi:e3b} and \eqref{tt3bi:e2}, by using the fact that  the centred random variables $(\tilde{\Delta} Z _{l,N}) ^{p} - \mathbf{E}  (\tilde{\Delta} Z _{l,N}) ^{p}$, $l\in\{1,..., \delta(N)\}$, are independent and identically distributed, one obtains that 
\begin{equation}
\label{tt3bi:e11}
{\rm Var}(\widetilde{v}_{N} ^{(1)})=\mathbf{E} (\widetilde{v}_{N} ^{(1)}) ^{2} =  \frac{2 ^{2pHN}\delta(N)}{2^{[N^\gamma]}}\mathbf{E}  \left[ (\tilde{\Delta} Z _{1,N}) ^{p} - \mathbf{E}  (\tilde{\Delta} Z _{1,N}) ^{p}\right]^{2}=2^{-[N^\nu]}\mathbf{E} (V_{N, p} ^{(1)}) ^{2},
\end{equation}
where the last equality follows from \eqref{e:add1}. Then combining \eqref{tt3bi:e10} and \eqref{tt3bi:e11} with \eqref{5i-1}, one gets that 
\begin{equation}
\label{tt3bi:e12}
\Big |{\rm Var}\big (\sqrt{m_{p}}\,W_{\frac{1}{2^{[N^\nu]}}}\big)-{\rm Var}\big (\widetilde{v}_{N} ^{(1)}\big)\Big |\le C 2 ^{-[N^\nu]+[N ^{\beta}]\frac{H-1}{q}}.
\end{equation}
On the other hand, by using \eqref{tt3bi:e3b} and \eqref{tt3bi:e2} and the same arguments which have allowed to obtain \eqref{20f-3}, it can be shown that 
\begin{equation}\label{tt3bi:e13}
	{\rm Var} \left(\langle D\widetilde{v}_{N} ^{(1)}, D(-L) ^{-1}\widetilde{v}_{N} ^{(1)}\rangle _{L ^{2}(\mathbb{R})}\right) \leq C 2 ^{-[N^\nu]-[N^{\gamma}]}.
\end{equation}
Then, it results from \eqref{tt3bi:e12}, \eqref{tt3bi:e13} and Theorem \ref{tt1} that 
\[
d_W \left (  \sqrt{m_{p}}\,W_{\frac{1}{2^{[N^\nu]}}}\,, \widetilde{v}_{N} ^{(1)} \right)\le C 2 ^{-[N^\nu]-[N^{\gamma}]},
\]
which, in view of \eqref{tt3bi:e9}, implies that 
\begin{equation}\label{tt3bi:e14}
d_W (\widetilde{B}_N, \sqrt{m_{p}}\,\mathcal{I}(h_N))\le C 2 ^{-[N^{\gamma}]}.
\end{equation}
Finally, since $0<\alpha <1$ and $0<\nu <\gamma\le \beta <1$, putting together \eqref{tt3bi:e30}, \eqref{tt3bi:e32}, \eqref{tt3bi:e35}, \eqref{tt3bi:e36} and \eqref{tt3bi:e14}, one obtains \eqref{tt3bi:e1}.
\qed

\begin{remark}
	A careful inspection of the proof of Theorem \ref{tt3bi} shows that, under the weaker condition that $h$ is a continuous function on $[0,1]$ (which may not be a H\"older function), one still has that $d_W(U_{N, p} (Z), \sqrt{m_{p}}\,\mathcal{I}(h))\le C \rho_h(2 ^{-N ^{\nu}})$, where the non-negative, non-decreasing continuous function $\rho_h$ is the uniform modulus of continuity of $h$, defined, for each $\delta\in [0,1]$, as $\rho_h (\delta)=\sup\big\{| h(x)-h(y)|,\, (x,y)\in [0,1]^2 \mbox{ and } |x-y|\le \delta\big\}$.
\end{remark}

\subsection{Volatility estimation}

Let $ h: [0,1] \to \mathbb{R}$  be a measurable function  satisfying (\ref{h}). We consider the stochastic process $ (X_{t}, t\in [0,1])$ defined as
\begin{equation}
	\label{x}
	X_{t}= x_{0}+ \int_{0} ^{t} h(s) dZ_{s},
\end{equation}
where $ x_{0} \in \mathbb{R}$ is deterministic and the stochastic integral with respect to $Z$ in (\ref{x}) is a Wiener-Hermite integral; recall that a brief presentation of such a stochastic integral has been given in Section \ref{sec21}. Observe that the function $h$ can be viewed as a measurable real-valued function defined on $\mathbb{R}$ which vanishes outside of the interval $[0,1]$ and satisfies, for every $t\in [0,1]$, $h \mathbbm{1}_{(0, t)} \in \vert \mathcal{H}\vert $; thus, one knows from Section \ref{sec21} that the Hermite-Wiener integral  $\int_{0} ^{t} h(s) dZ_{s}$ is well-defined. One mentions in passing that the fact that, for all $t\in [0,1]$, $h \mathbbm{1}_{(0, t)} \in \vert \mathcal{H}\vert $ is a straightforward consequence of the boundedness of the function $h$ (see(\ref{h})). Indeed, since $\|h\|_\infty<\infty$ and $H\in (1/2,1)$, one has, for each $t\in [0,1]$,
\begin{equation*}
	\int_{0} ^{t} \int_{0} ^{t} \vert h(u)\vert \cdot \vert h(v)\vert \cdot \vert u-v \vert ^{2H-2} dudv \leq \|h\|_\infty^2 \int_{0} ^{1} \int_{0} ^{1} \vert u-v \vert ^{2H-2} dudv<\infty.
\end{equation*} 

Let us emphasize that the estimation of the integrated volatility $\int_{0} ^{1} h(s) ^{p} ds$, with $p\geq 1$, constitutes a common topic in financial mathematics. When the noise $Z$ in (\ref{x}) is a semimartingale, there exists a huge literature on this topic. We refer, among many others, to \cite{BS}, \cite{Betal}, \cite{Ja1} or  \cite{PV} for several approaches to estimate the integrated volatility, even when it is assumed to be random. When the process $X$ is observed at discrete times on the interval $[0, 1]$, one of the most usual methods to do it  is via a power variation of the process $X$.  We will use a similar approach in the case of the Hermite noise in (\ref{x}), but our power variation will be defined by using only some special increments of the process $X$ given by (\ref{x}). This will allow to prove that the power variation sequence defined below constitutes a consistent estimator for the integrated volatility and moreover, it satisfies a Central Limit Theorem after a proper renormalization. Related results obtained for the case when $Z$ is the fractional Brownian motion can be found in \cite{BCP} or \cite{CNW}.

We start by defining the modified power variation of the process $X$.  We set, for any integer numbers $p\geq 2$ and $N\geq 1$,
\begin{equation}
	\label{snpx}
	S_{N, p}(X) = \frac{ 2 ^{pHN}}{2 ^{[N ^{\gamma}]}} \sum _{l=1} ^{2 ^{ [N ^{\gamma}]}}(\Delta X_{l, N}) ^{p}, 
\end{equation}
where, for $ l=1,..., 2 ^{ [N ^{\gamma}]}$, 
\begin{equation}\label{incrementsforsnpx}
	\Delta X_{l, N}= X _{ \frac{ l}{ 2 ^{ [N ^{\gamma}]}}+ 2 ^{-N} }- X _{ \frac{ l}{ 2 ^{ [N ^{\gamma}]}}}.
\end{equation}
Thus, similarly to what has been done \eqref{dz} and \eqref{snpz}, we use the increments of $X$ of length $2 ^{-N}$  located at the points $\frac{l} {2^{[N ^{\gamma}]}}$ of the unit interval $[0,1]$. The following result shows that the modified power variation of $X$ provides a consistent estimator for the integrated volatility.

\begin{prop}
	Assume $h$ satisfies (\ref{h}) and consider the sequence $(S_{N, p}(X), N\geq 1)$ given by (\ref{snpx})  with an integer $p\geq 2$. Then, setting $\mu_{p}= \mathbf{E} Z_{1} ^{p}$, one has
	\begin{equation}\label{eqn:convergencelimite}
		S_{N, p}(X) \to _{N \to \infty} \mu _{p} \int_{0} ^{ 1} h(s) ^{p} ds  \mbox{ in } L ^{1}(\Omega).
	\end{equation}
\end{prop}
{\bf Proof: } As we already did it in the proof of Theorem \ref{tt3bi}, we use an approach based on the approximation of the stochastic integral in (\ref{x}) by Riemann sums  inspired by the works \cite{BCP}  or \cite{CNW}. Let $M$ and $N$ be two arbitrary integer numbers such that $2<2^{(1-\beta)^{-1}}\le M \leq [N ^{\gamma}]$.  We decompose the difference $S_{N, p}(X) - \mu _{p} \int_{0} ^{ 1} h(s) ^{p} ds$ into four  terms as follows:
\begin{eqnarray}
&&	S_{N, p}(X) -\mu _{p} \int_{0} ^{ 1} h(s) ^{p} ds\nonumber\\
&&= \frac{ 2 ^{pHN}}{2 ^{[N ^{\gamma}]}} \sum _{l=1} ^{2 ^{ [N ^{\gamma}]}}(\Delta X_{l, N}) ^{p}-  \frac{ 2 ^{pHN}}{2 ^{[N ^{\gamma}]}} \sum _{l=1} ^{2 ^{ [N ^{\gamma}]}} h\left( \frac{l}{ 2 ^{N ^{\gamma}]}}\right) ^{p} (\Delta Z_{l, N}) ^{p}\nonumber \\
&&+  \frac{ 2 ^{pHN}}{2 ^{[N ^{\gamma}]}} \sum _{k=1} ^{2 ^{M}} \sum _{l=(k-1) \frac{ 2 ^{[ N^{\gamma}]}}{2 ^{M}}+1}^{l= k  \frac{ 2 ^{[ N ^{\gamma}]}}{2 ^{M}}}\left( h\left( \frac{l}{ 2 ^{[N ^{\gamma}]}}\right) ^{p}- h\left( \frac{k-1}{2 ^{M}}\right) ^{p} \right) (\Delta Z_{l, N}) ^{p} \nonumber\\
&&+   \frac{ 2 ^{pHN}}{2 ^{[N ^{\gamma}]}} \sum _{k=1} ^{2 ^{M}}h\left( \frac{k-1}{2 ^{M}}\right) ^{p} \sum _{l=(k-1) \frac{ 2 ^{[ N^{\gamma}]}}{2 ^{M}}+1}^{l= k  \frac{ 2 ^{[ N ^{\gamma}]}}{2 ^{M}}}(\Delta Z_{l, N}) ^{p} -\mu _{p} \frac{1}{ 2 ^{M}}\sum _{k=1} ^{2 ^{M}}h\left( \frac{k-1}{2 ^{M}}\right) ^{p}\nonumber\\
&&+ \mu _{p} \left(  \frac{1}{ 2 ^{M}}\sum _{k=1} ^{2 ^{M}}h\left( \frac{k-1}{2 ^{M}}\right) ^{p}-\int_{0} ^{1}h(s)^{p} ds \right)\nonumber\\
&:=& A_{N, p}+ B_{N, M, p}+ C_{N, M, p} + D_{M, p}.\label{14i-1}
\end{eqnarray}
Our next goal is to estimate from above each one of the latter four summands.
\vskip0.2cm

\noindent{\it Estimation of $ A_{N, p}$. }We use the triangle inequality and the inequality
\begin{equation}\label{10i-2}
	\vert A ^{p}- B ^{p} \vert \leq 2 ^{p-2}\, p \left( |B| ^{p-1} \vert A-B\vert + \vert A-B\vert ^{p} \right),
\end{equation}
which holds, for any $A, B\in \mathbb{R}$ and for all integer $p\geq 2$. Thus, we  obtain 
\begin{equation*}
	\vert A_{N, p}\vert \leq C_{p} \left( \vert A_{N, p, 1}\vert + \vert A_{N, p, 2}\vert\right),
\end{equation*}
with

\begin{equation*}
	A_{N, p, 1}=  \frac{ 2 ^{pHN}}{2 ^{[N ^{\gamma}]}} \sum _{l=1} ^{2 ^{ [N ^{\gamma}]}}\left| h\left( \frac{l}{ 2 ^{[N ^{\gamma}]}}\right)\right|  ^{p-1} \vert \Delta Z_{l, N}\vert ^{p-1} \left |\Delta X_{l, N} - h\left( \frac{l}{ 2 ^{[N ^{\gamma}]}}\right)\Delta Z_{l, N}\right |
\end{equation*}
and
\begin{equation*}
	A_{N, p,2}= \frac{ 2 ^{pHN}}{2 ^{[N ^{\gamma}]}} \sum _{l=1} ^{2 ^{ [N ^{\gamma}]}}\left| \Delta X_{l, N} - h\left( \frac{l}{ 2 ^{[N ^{\gamma}]}}\right)\Delta Z_{l, N} \right| ^{p}.
\end{equation*}

Let $a\geq 2$ be an arbitrary positive integer number, it follows from \eqref{x} that 
\begin{eqnarray*}
\mathbf{E} 	\left| \Delta X_{l, N} - h\left( \frac{l}{ 2 ^{[N ^{\gamma}]}}\right)\Delta Z_{l, N} \right| ^{a}
	= \mathbf{E}  \left| \int_{ \frac{ l}{ 2 ^{[ N ^{\gamma}]}}}^{ \frac{ l}{ 2 ^{[ N ^{\gamma}]}}+ 2^{ -N}}\left( h (s)- h\left( \frac{l}{ 2 ^{[N ^{\gamma}]}}\right)\right) dZ_{s} \right| ^{a}.
\end{eqnarray*}
Since the Wiener integral of any function of $\vert \mathcal{H}\vert$ with respect to $Z$ is an element of the $q$th Wiener chaos, using the hypercontractivity property (\ref{hyper})  and the isometry of the Wiener integral (\ref{iso-w}), one can derive from the previous equality that

	\begin{eqnarray*}
		&&\mathbf{E} 	\left| \Delta X_{l, N} - h\left( \frac{l}{ 2 ^{[N ^{\gamma}]}}\right)\Delta Z_{l, N} \right| ^{a}\\
		&&\leq C \left( \mathbf{E}  \left| \int_{ \frac{ l}{ 2 ^{[ N ^{\gamma}]}}}^{ \frac{ l}{ 2 ^{[ N ^{\gamma}]}}+ 2^{ -N}}\left( h(s)- h\left( \frac{l}{ 2 ^{[N ^{\gamma}]}}\right)\right) dZ_{s} \right| ^{2}\right) ^{\frac{a}{2}}\\
		&=& \left( C \int_{ \frac{ l}{ 2 ^{[ N ^{\gamma}]}}}^{ \frac{ l}{ 2 ^{[ N ^{\gamma}]}}+ 2^{ -N}}\int_{ \frac{ l}{ 2 ^{[ N ^{\gamma}]}}}^{ \frac{ l}{ 2 ^{[ N ^{\gamma}]}}+ 2^{ -N}} 
		\left( h(u)- h\left( \frac{l}{ 2 ^{[N ^{\gamma}]}}\right)\right)\right.\\
		&&\left. \times  \left( h(v)- h\left( \frac{l}{ 2 ^{[N ^{\gamma}]}}\right)\right)  \vert u-v\vert ^{2H-2} dudv\right) ^{\frac{a}{2}}\\
		&\leq & C\left( 2 ^{-2\alpha N}  \int_{ \frac{ l}{ 2 ^{[ N ^{\gamma}]}}}^{ \frac{ l}{ 2 ^{[ N ^{\gamma}]}}+ 2^{ -N}} \int_{ \frac{ l}{ 2 ^{[ N ^{\gamma}]}}}^{ \frac{ l}{ 2 ^{[ N ^{\gamma}]}}+ 2^{ -N}}\vert u-v\vert ^{2H-2} dudv \right) ^{\frac{a}{2}},
\end{eqnarray*}
where we used (\ref{h}) for the last inequality. Consequently, we obtain 
\begin{equation}\label{10i-1}
	\mathbf{E} 	\left| \Delta X_{l, N} - h\left( \frac{l}{ 2 ^{[N ^{\gamma}]}}\right)\Delta Z_{l, N} \right| ^{a}\leq C 2 ^{-\alpha Na}2 ^{-HNa}.
\end{equation}
We will use Cauchy-Schwarz's inequality and (\ref{10i-1}) to bound the terms denoted by $A_{N, p, 1}$ and $ A_{N, p, 2}$. For $ A_{ N, p, 1}$, we write 
\begin{eqnarray*}
	\mathbf{E}  \vert A_{N, p, 1}\vert &\leq & C \frac{ 2 ^{pHN}}{2 ^{[N ^{\gamma}]}} \sum _{l=1} ^{2 ^{ [N ^{\gamma}]}}\mathbf{E}  \vert \Delta Z_{l, N}\vert ^{p-1} \left| \Delta X_{l, N} - h\left( \frac{l}{ 2 ^{[N ^{\gamma}]}}\right)\Delta Z_{l, N}\right| \\
	&\leq & C  \frac{ 2 ^{pHN}}{2 ^{[N ^{\gamma}]}} \sum _{l=1} ^{2 ^{ [N ^{\gamma}]}}\left( \mathbf{E}  \vert \Delta Z_{l, N}\vert ^{2p-2}\right) ^{\frac{1}{2}}  \left( \mathbf{E} \left| \Delta X_{l, N} - h\left( \frac{l}{ 2 ^{[N ^{\gamma}]}}\right)\Delta Z_{l, N}\right|^2 \right) ^{\frac{1}{2}}\\
	&\leq & C   \frac{ 2 ^{pHN}}{2 ^{[N ^{\gamma}]}} \sum _{l=1} ^{2 ^{ [N ^{\gamma}]}}2 ^{-(p-1)HN} 2 ^{-HN} 2 ^{-\alpha N}\\
	&\leq & C 2 ^{-\alpha N}.
\end{eqnarray*}
Regarding  $ A_{N, p, 2}$, 
\begin{eqnarray*}
	\mathbf{E}  \vert A_{N, p, 2}\vert 
	&\leq & \frac{ 2 ^{pHN}}{2 ^{[N ^{\gamma}]}} \sum _{l=1} ^{2 ^{ [N ^{\gamma}]}}\left( \mathbf{E}  \left| \Delta X_{l, N} - h\left( \frac{l}{ 2 ^{[N ^{\gamma}]}}\right)\Delta Z_{l, N} \right| ^{2p}\right) ^{\frac{1}{2}}\\
	&\leq &C \frac{ 2 ^{pHN}}{2 ^{[N ^{\gamma}]}} \sum _{l=1} ^{2 ^{ [N ^{\gamma}]}} 2 ^{ -\alpha pN} 2 ^{-HpN}\leq C 2 ^{-\alpha pN}.
\end{eqnarray*}
We then conclude  that 
\begin{equation}\label{14i-2}
	\mathbf{E} \vert A_{N, p} \vert\leq C 2 ^{-\alpha N}.
\end{equation}

\vskip0.2cm
\noindent {\it Estimation of $B_{N, M, p}$. } We notice that for every integer $l\in \left[ (k-1) \frac{ 2 ^{[N ^{\gamma}]}}{2 ^{M}}+1, k \frac{ 2 ^{[N ^{\gamma}]}}{2 ^{M}}\right]$, 
\begin{equation}\label{13i-1}
	\left| \frac{l}{2 ^{[N ^{\gamma}]}}- \frac{k-1}{2 ^{M}}\right|\leq 2 ^{-M}.
\end{equation}
Then, by (\ref{h}) and the inequality (\ref{10i-2}), 
\begin{equation*}
\left| 	h\left( \frac{l}{ 2 ^{[N ^{\gamma}]}}\right) ^{p}- h\left( \frac{k-1}{2 ^{M}}\right) ^{p}\right|\leq C 2 ^{-M\alpha}, 
\end{equation*}
for any integer number $p\geq 2$. It gives
\begin{eqnarray}
\label{e:b-B}
\mathbf{E} \vert B_{N, M, p}\vert & \leq& 	\frac{ 2 ^{pHN}}{2 ^{[N ^{\gamma}]}} \sum _{k=1} ^{2 ^{M}} \sum _{l=(k-1) \frac{ 2 ^{[ N ^{\gamma}]}}{2 ^{M}}+1}^{l= k  \frac{ 2 ^{[ N ^{\gamma}]}}{2 ^{M}}}\left| h\left( \frac{l}{ 2 ^{[N ^{\gamma}]}}\right) ^{p}- h\left( \frac{k-1}{2 ^{M}}\right) ^{p} \right|  \mathbf{E} \vert \Delta Z_{l, N}\vert  ^{p}\nonumber \\
&\leq & C 	\frac{ 2 ^{pHN}}{2 ^{[N ^{\gamma}]}} \sum _{k=1} ^{2 ^{M}} \sum _{l=(k-1) \frac{ 2 ^{[N^{\gamma}]}}{2 ^{M}}+1}^{l= k  \frac{ 2 ^{[ N ^{\gamma}]}}{2 ^{M}}} 2 ^{-M \alpha }  \mathbf{E} \vert \Delta Z_{l, N}\vert  ^{p}\nonumber\\
&\leq & C 	\frac{ 2 ^{pHN}}{2 ^{[N ^{\gamma}]}} \sum _{k=1} ^{2 ^{M}} \sum _{l=(k-1) \frac{ 2 ^{[ N ^{\gamma}]}}{2 ^{M}}+1}^{l= k  \frac{ 2 ^{[ N ^{\gamma}]}}{2 ^{M}}}  2 ^{-M \alpha } 2 ^{-pHN}\nonumber\\
&\leq & C 2 ^{ -M \alpha }. 
\end{eqnarray}

\vskip0.2cm
\noindent {\it Estimation of $ C_{N, M, p}$. }Using  the equalities $\mu_p=\mathbf{E}Z_1^p= 2 ^{pHN}\mathbf{E}\left(\Delta Z_{l, N}) ^{p} \right)$, for all $l\in\{1,\ldots, 2^{[N^\gamma]}\}$, we can express this summand as
\begin{eqnarray}
\label{e1:b-C}
	C_{N, M, p} &=&\frac{1}{ 2 ^{M}}  \sum _{k=1} ^{2 ^{M}}h\left( \frac{k-1}{2 ^{M}}\right) ^{p} 	\frac{ 2 ^{pHN}2 ^{M}}{2 ^{[N ^{\gamma}]}}\sum _{l=(k-1) \frac{ 2 ^{[ N ^{\gamma}]}}{2 ^{M}}+1}^{l= k  \frac{ 2 ^{[ N ^{\gamma}]}}{2 ^{M}}} \left( (\Delta Z_{l, N})^{p}- \mathbf{E}  (\Delta Z_{l, N}) ^{p} \right)\nonumber\\
	&=& C_{N, M, p} ^{(1)} + R_{N, M, p},
\end{eqnarray}
with
\begin{equation*}
	C_{N, M, p} ^{(1)} =\frac{1}{ 2 ^{M}}  \sum _{k=1} ^{2 ^{M}}h\left( \frac{k-1}{2 ^{M}}\right) ^{p} 	\frac{ 2 ^{pHN}2 ^{M}}{2 ^{[N ^{\gamma}]}}\sum _{l=(k-1) \frac{ 2 ^{[ N ^{\gamma}]}}{2 ^{M}}+1}^{l= k  \frac{ 2 ^{[ N ^{\gamma}]}}{2 ^{M}}} \left( (\tilde{\Delta} Z_{l, N})^{p}- \mathbf{E}  (\tilde{\Delta} Z_{l, N})\right)
\end{equation*}
and

\begin{eqnarray*}
	R_{N, M, p}= &&\frac{1}{ 2 ^{M}}  \sum _{k=1} ^{2 ^{M}}h\left( \frac{k-1}{2 ^{M}}\right) ^{p} 	\frac{ 2 ^{pHN}2 ^{M}}{2 ^{[N ^{\gamma}]}}\\
	&&\hspace{0.3cm}\times\sum _{l=(k-1) \frac{ 2 ^{[ N ^{\gamma}]}}{2 ^{M}}+1}^{l= k  \frac{ 2 ^{[ N ^{\gamma}]}}{2 ^{M}}} \sum_{j=0} ^{p-1} \binom{p}{j}  \left ( (\tilde{\Delta } Z_{l, N})^{j} (\check{\Delta}Z_{l, N})^{p-j}-\mathbf{E}(\tilde{\Delta } Z_{l, N})^{j} (\check{\Delta}Z_{l, N})^{p-j} \right).
\end{eqnarray*}
Since $h$ is a bounded function on the interval $[0,1]$ (see (\ref{h})), similarly to the proof of (\ref{8i-2}), it can be shown that
\begin{equation}
\label{e2:b-C}
	\mathbf{E}  \vert R_{N, M, p}\vert \leq C 2 ^{[ N ^{\beta}]\frac{H-1}{q}},
\end{equation}
for all integers $N$ and $M$ with $2^{(1-\beta)^{-1}}\le M\leq [N ^{\gamma}]$. We next deal with $ C_{N, M, p} ^{(1)}$. We will compute its norm in $ L^{2}(\Omega)$, in order to benefit from the independence property in Proposition \ref{pp1}. We have, for all integers $M, N$ with $2^{(1-\beta)^{-1}}\le M\leq [N ^{\gamma}]$, 
\begin{eqnarray*}
&&	\mathbf{E}  ( C_{N, M, p}^{(1)})^{2} \\
&=&\frac{1}{ 2 ^ {2M}}  \sum _{k=1} ^{2 ^{M}}h\left( \frac{k-1}{2 ^{M}}\right) ^{p} 	\frac{ 2 ^{2pHN}2 ^{2M}}{2 ^{2[N ^{\gamma}]}}\sum _{l=(k-1) \frac{ 2 ^{[ N ^{\gamma}]}}{2 ^{M}}+1}^{l= k  \frac{ 2 ^{[ N ^{\gamma}]}}{2 ^{M}}}\mathbf{E}  \left( (\tilde{\Delta} Z_{l, N})^{p}- \mathbf{E}  (\tilde{\Delta} Z_{l, N})^p\right)^{2}\\
	& \leq&  C \frac{ 2 ^{2pHN}}{2 ^{2[N ^{\gamma}]}}\sum _{k=1} ^{2 ^{M}}\sum _{l=(k-1) \frac{ 2 ^{[ N ^{\gamma}]}}{2 ^{M}}+1}^{l= k  \frac{ 2 ^{[ N ^{\gamma}]}}{2 ^{M}}}\mathbf{E}  \left( (\tilde{\Delta} Z_{l, N})^{p}- \mathbf{E}  (\tilde{\Delta} Z_{l, N})^p\right)^{2}\\
	&\leq & C \frac{ 1}{2 ^{2[N ^{\gamma}]}}\sum _{k=1} ^{2 ^{M}}\sum _{l=(k-1) \frac{ 2 ^{[ N ^{\gamma}]}}{2 ^{M}}+1}^{l= k  \frac{ 2 ^{[ N ^{\gamma}]}}{2 ^{M}}} 1\leq  C 2 ^{- [N ^{\gamma}]}.
\end{eqnarray*}
Thus, using the latter inequality, Cauchy-Schwarz's inequality, \eqref{e1:b-C} and \eqref{e2:b-C} we get that
\begin{equation}
\label{e:b-C}
	\mathbf{E}  \vert C_{N, M, p} \vert \leq C 2 ^{-\frac{ [ N ^{\gamma}]}{2}}.
\end{equation}

\vskip0.2cm

\noindent {\it Estimation of $ D_{M, p}$. } We can write 
\begin{equation*}
	D_{M, p}= \mu _{p} \sum _{k=1} ^{2^{M}} \int_{\frac{k-1}{2^{M}}} ^{\frac{k}{2^{M}}}\left( h\left( \frac{k-1}{2^{M}}\right) ^{p} -h(s) ^{p} \right) ds,  
\end{equation*}
and by (\ref{10i-2}) and (\ref{h}), 
\begin{eqnarray}
\label{e:b-D}
	\mathbf{E}  \vert D_{M, p}\vert &\leq& C \sum _{k=1} ^{2^{M}} \int_{\frac{k-1}{2^{M}}} ^{\frac{k}{2^{M}}}\left| h\left( \frac{k-1}{2^{M}}\right) ^{p} -h(s) ^{p} \right| ds\nonumber\\
	&\leq &C \sum _{k=1} ^{2^{M}} \int_{\frac{k-1}{2^{M}}} ^{\frac{k}{2^{M}}} 2 ^{-M\alpha} ds \leq C 2 ^{-M\alpha}.
\end{eqnarray}
Finally, using \eqref{14i-1}, \eqref{14i-2}, \eqref{e:b-B}, \eqref{e:b-C} and \eqref{e:b-D}, we obtain, for all integer $M\ge 2^{(1-\beta)^{-1}}$, 
\[
\limsup_{N\rightarrow +\infty}\mathbf{E}\left |S_{N, p}(X) -\mu _{p} \int_{0} ^{ 1} h(s) ^{p} ds\right| \leq C 2 ^{-M\alpha}.
\]
Thus, letting $M$ goes to $+\infty$, it follows that 
\[
\lim_{N\rightarrow +\infty}\mathbf{E}\left |S_{N, p}(X) -\mu _{p} \int_{0} ^{ 1} h(s) ^{p} ds\right| =\limsup_{N\rightarrow +\infty}\mathbf{E}\left |S_{N, p}(X) -\mu _{p} \int_{0} ^{ 1} h(s) ^{p} ds\right| =0.
\]

\qed

For every integers $p\geq 2$ and $N\ge 1$, we set
\begin{equation}
	\label{vnpx}
	V_{N, p}(X)= \sqrt{ 2 ^ { [N ^{\gamma}]}} \left( S_{N, p}(X)- \mu_{p} \int_{0} ^{1} h(s) ^{p} ds\right), 
\end{equation}
where $S_{N, p}(X)$ is given by (\ref{snpx}). We recall that $h: [0,1] \to \mathbb{R}$ is a deterministic function satisfying (\ref{h}). The next result shows that, as soon as the function $h$ is smooth enough, the renormalized modified power variation $V_{N, p}(X)$ converges in distribution to a Gaussian random variable which can be expressed, up to the multiplicative constant $\sqrt{m_{p}}$, as a Wiener integral of the function $h$.

\begin{theorem}
Let $p\geq 2$ be an integer number. 	Consider the sequence $(V_{N, p}(X), N\geq 1)$ given by (\ref{vnpx}) and assume that $h$ satisfies (\ref{h}) with $ \alpha > \frac{1}{2}$. Then 
	\begin{equation*}
		V_{N, p} (X) \to ^{(d)} _{N \to \infty} \sqrt{m_{p}}\int_{0} ^{1} h(s) dW_{s},
	\end{equation*}
where $ (W_{t}, t\geq 0)$ is a Wiener process and $m_{p}$ is given by (\ref{mp}). Moreover, for each $\nu\in (0,\gamma)$, there is a constant $C$, which depends on $\nu$ and $h$, such that, for all $N$ large enough, one has 
\begin{equation}
	\label{tt3bi:e1b}
	d_W\left(V_{N, p} (Z), \sqrt{m_{p}}\int_{0} ^{1} h(s)^p dW_{s}\right)\le C 2 ^{-\alpha N ^{\nu}}.
\end{equation}
\end{theorem}
{\bf Proof: }We write
\begin{eqnarray*}
	V_{N, p}(X)&=&  \sqrt{ 2 ^ { [N ^{\gamma}]}}\left( \frac{ 2 ^{pHN}}{ 2 ^{[N ^{\gamma}]}} \sum _{l=1} ^{2 ^{ [N ^{\gamma}]}}(\Delta X_{l, N}) ^{p} -\mu_{p} \int_{0} ^{1}h(s) ^{p} ds  \right) \\
	&=&  \sqrt{ 2 ^ { [N ^{\gamma}]}} \left(\frac{ 2 ^{pHN}}{ 2 ^{[N ^{\gamma}]}} \sum _{l=1} ^{2 ^{ [N ^{\gamma}]}}(\Delta X_{l, N}) ^{p}-  \frac{ 2 ^{pHN}}{ 2 ^{[N ^{\gamma}]}} \sum _{l=1} ^{2 ^{ [N ^{\gamma}]}}h\left( \frac{ l}{2 ^{ [N ^{\gamma}]}} \right)^{p} (\Delta Z_{l, N}) ^{p} \right)\\
	&&+   \sqrt{ 2 ^ { [N ^{\gamma}]}} \left(  \frac{ 2 ^{pHN}}{ 2 ^{[N ^{\gamma}]}} \sum _{l=1} ^{2 ^{ [N ^{\gamma}]}}h\left( \frac{ l}{2 ^{ [N ^{\gamma}]}} \right)^p (\Delta Z_{l, N}) ^{p} - \frac{1}{ 2 ^{[N ^{\gamma}]}}\mu_{p}  \sum _{l=1} ^{2 ^{ [N ^{\gamma}]}}h\left( \frac{ l}{2 ^{ [N ^{\gamma}]}} \right)^{p} \right)\\
	&&+   \sqrt{ 2 ^ { [N ^{\gamma}]}} \left(  \frac{1}{ 2 ^{[N ^{\gamma}]}}\mu_{p}  \sum _{l=1} ^{2 ^{ [N ^{\gamma}]}}h\left( \frac{ l}{2 ^{ [N ^{\gamma}]}} \right)^{p}-\mu_{p} \int_{0} ^{1} h(s) ^{p} ds \right)\\
	&:=& T_{1, N, p}+ T_{2, N, p} + T_{3, N, p}.
\end{eqnarray*}
Let us estimate each of the three summands from above. First, we notice that 
\begin{equation*}
	T_{1, N, p}=   \sqrt{ 2 ^ { [N ^{\gamma}]}}A_{ N, p},
\end{equation*}
where $ A_{N, p}$ is given by (\ref{14i-1}). By using the bound (\ref{14i-2}), we obtain

\begin{equation}\label{9d-1}
	\mathbf{E}  \vert T_{1, N, p} \vert \leq C 2 ^{\frac{ [N ^{\gamma}]}{2}-\alpha N} \to_{N \to \infty } 0. 
\end{equation}

For $ T_{2, N, p}$, we write
\begin{eqnarray*}
	T_{2, N, p} &=& = \sqrt{ 2 ^ { [N ^{\gamma}]}}   \frac{ 2 ^{pHN}}{ 2 ^{[N ^{\gamma}]}}  \sum _{l=1} ^{2 ^{ [N ^{\gamma}]}}h \left( \frac{ l}{2 ^{ [N ^{\gamma}]}} \right)^{p} \left( (\Delta Z_{l, N} ) ^{p} -2 ^{-pHN} \mu_{p}\right) \\
	&=&   \frac{ 2 ^{pHN}}{ \sqrt{2 ^{[N ^{\gamma}]}}  } \sum _{l=1} ^{2 ^{ [N ^{\gamma}]}}h\left( \frac{ l}{2 ^{ [N ^{\gamma}]}} \right)^{p}  \left( (\Delta Z_{l, N} ) ^{p} -\mathbf{E}  (\Delta Z_{l, N}) ^{p} \right).
\end{eqnarray*}
Since $h$ is $\alpha$-H\"older continuous, $ h^{p}$ is also $\alpha$-H\"older continuous on $[0, 1]$. Thus, we can apply Theorem \ref{tt3bi} to conclude that
\begin{equation*}
	T_{2, N, p}\to ^{(d)} _{N \to \infty} \sqrt{m_{p}} \int_{0} ^{1} h(s) ^{p}dW_{s},
\end{equation*}
where $(W_{t}, t\geq 0)$ denotes a Brownian motion. Moreover, for $N $ large enough,
\begin{equation}\label{9d-2}
	d_{W} \left( T_{2, N, p}, \sqrt{m_{p}} \int_{0} ^{1} h(s) ^{p}dW_{s}\right)\leq C 2 ^{-\alpha N ^{\nu}}.
\end{equation}

 We finally treat the summand $ T_{3, N, p}$. This term can be estimated as follows: 
\begin{eqnarray*}
	T_{3, N, p}= \sqrt{ 2 ^ { [N ^{\gamma}]}}  \mu_{p} \sum _{l=1} ^{2 ^{ [N ^{\gamma}]}} 
	\int_{ \frac{ l-1} {2 ^{[N ^{\gamma}]}}} ^{ \frac{ l} {2 ^{[N ^{\gamma}]}}}\left( h\left( \frac{ l}{2 ^{ [N ^{\gamma}]}} \right)^{p} -h(s) ^{p} ds \right) ds 
\end{eqnarray*}
and then, via(\ref{h}), 
\begin{eqnarray}
	\mathbf{E} \vert T_{3, N, p}\vert &\leq & \sqrt{ 2 ^ { [N ^{\gamma}]}}  \mu_{p} \sum _{l=1} ^{2 ^{ [N ^{\gamma}]}} 
	\int_{ \frac{ l-1} {2 ^{[N ^{\gamma}]}}} ^{ \frac{ l} {2 ^{[N ^{\gamma}]}}}\left| h\left( \frac{ l}{2 ^{ [N ^{\gamma}]}} \right)^{p} -h(s) ^{p} ds \right| ds \nonumber\\
	&\leq & C \sqrt{ 2 ^ { [N ^{\gamma}]}} 2 ^{-\alpha [ N^{\gamma}]}\to _{N \to \infty} 0,\label{9d-3}
\end{eqnarray}
since $\alpha >\frac{1}{2}$. 

The desired bound for the Wasserstein distance follows from (\ref{9d-1}), (\ref{9d-2}) and (\ref{9d-3}). \qed 

\subsection{Numerical experiments}

The goal of this last subsection is to make a numerical study of the performances of the estimator of integrated volatility introduced in \eqref{snpx} by using data issued from simulated sample paths of the stochastic process $X=(X_t, t\in [0,1])$ in \eqref{x}. For the sake of simplicity one assumes that $x_0=0$. Thus, for all $t\in [0,1]$, the random variable $X_{t}$, which is sometimes denoted by $X(t)$, reduces to
\begin{equation}
	\label{x:bis}
	X_{t}= X(t)=\int_{0} ^{t} h(s) dZ_{s}.
\end{equation}
\begin{remark}
\label{rem:you}
For numerical approximation purposes, one needs that, for each $t\in [0,1]$, the Wiener-Hermite integral $\int_0^t h(s) dZ_s$ (see Section \ref{sec21}) be also a well-defined pathwise Riemann-Stieltjes integral. To this end, one imposes to the deterministic function $h$ to satisfy the following assumption:
\begin{itemize}
\item[($\mathcal{A}$)] The function $h$ is $\alpha$-H\"older continuous on the interval $[0,1]$ of some order $\alpha\in (0,1]$ such that 
\begin{equation}
\label{con-al}
\alpha+H>1,
\end{equation}
\end{itemize}
where $H$ denotes the Hurst parameter of the Hermite process $Z$. Actually, when $\delta<H$ is close enough to $H$, \eqref{con-al} implies that $\alpha+\delta>1$ then, since sample paths of the Hermite process $Z=(Z_s, s\in [0,1])$ are $\delta$-H\"older continuous functions, one knows from e.g. \cite{LCL04} that, for each $t\in [0,1]$, the pathwise Riemann-Stieltjes integral on $[0,t]$ of $h$ with respect to the integrator $Z$ is well-defined. Moreover,  in virtue of \cite[Proposition 3.2]{T3}, one knows that the latter integral coincides with the Wiener-Hermite integral $\int_0^t h(s) dZ_s$. Thus, for any $t \in [0,1]$ and any sequence
\[ \left(P_j:=\left\{ 0 =x_0^{(j)}<x_1^{(j)}< \dots < x_m^{(j)}=t\right\}\right)_j\]
of subdivisions of interval $[0,t]$ whose mesh tends to $0$, one has, almost surely,
\begin{equation}\label{eqn:riemstiint}
\int_0^t h(s) dZ_s = \lim_{j \to + \infty} \sum_{k=0}^{m-1} h(x_k^{(j)}) \left(Z_{x_{k+1}^{(j)}}-Z_{x_{k}^{(j)}} \right).
\end{equation}
\end{remark}


In view of \eqref{x:bis} and \eqref{eqn:riemstiint}, loosely speaking, a reasonable strategy to simulate a sample path of the process $X=(X_t, t\in [0,1])$ would consist in considering a partition of $[0,1]$ and computing from it the Riemann-Stieltjes sum appearing in the right-hand side of \eqref{eqn:riemstiint} thanks to a simulated sample path of the Hermite process $Z$ itself, obtained by using the simulation method which was introduced very recently in \cite{AHL2}.

Let us now concisely describe the latter simulation method.
It relies on the resolution parameter $J \in \N$, which corresponds to the level of approximation (see equation \eqref{rateofconvergence} below). Also, it relies on the two other parameters $a \in (1/2,1)$ and $\varepsilon>0$, which are needed to insure the rate of convergence \eqref{rateofconvergence} for the approximation procedure, see Section 3 in \cite{AHL2} and references to \cite{AHL1} therein. Having introduced these three parameters, let us fix some notations. %
We recall that $q\in\N$ is the order of the Wiener chaos to which the Hermite process $Z$ of Hurst parameter $H\in (1/2,1)$ belongs. In the sequel we set $\delta=\frac{H-1}{q}+\frac{1}{2}$. 
The univariate Meyer fractional scaling function \cite{Meyer1999} of order $\delta$ is denoted by  $\Phi_\Delta^{(\delta)}$ and defined through its Fourier transform as
\[ \widehat{\Phi}_\Delta^{(\delta)}(\xi) = \left( \frac{1-e^{-i \xi}}{i \xi} \right)^\delta \widehat{\phi}(\xi)~\forall \, \xi\in\R\setminus\{0\} \text{ and }  \widehat{\Phi}_\Delta^{(\delta)}(0)=1,\]
where $\phi$ is a univariate Meyer scaling function \cite{LemarieRieusset1986}. Let $(g_{J,k}^\phi)_{ k \in \Z}$ be the sequence of i.i.d. $\mathcal{N}(0,1)$ Gaussian random variable defined, for all $k \in \Z$, as
\[ g_{J,k}^\phi := 2^{J/2} I_1 \left(\phi(2^J \bullet - k) \right),\]
the Gaussian FARIMA $(0,\delta,0)$ sequence $(Z_{J,\ell}^{(\delta)})_{\ell \in \Z}$ associated to $(g_{J,k}^\phi)_{ k \in \Z}$ is given, for all $\ell \in \Z$, by
\[ Z_{J,\ell}^{(\delta)}:= g_{J,\ell}^\phi  + \sum_{p=1}^{+ \infty} \gamma_p^{(\delta)} g_{J,\ell-p}^\phi, \]
where the coefficients $\gamma_p^{(\delta)} = \frac{\delta \, \Gamma(p+\delta)}{\Gamma(p+1) \Gamma(\delta+1)}$, for all $p \in \N$. We mention in passing that $(Z_{J,\ell}^{(\delta)})_{\ell \in \Z}$ can be simulated by using, for instance, a circulant matrix embedding procedure, which is an exact method \cite{daviesharte,woodchan}. 

Our next goal is to introduce a piecewise linear stochastic process which approximates the Hermite process $Z$. To this end, let us now consider the infinite set
\[
\mathcal{I}_J := \N \cap [2 ^{J(1-a)},+ \infty),
\]
and, for any $m \in \mathcal{I}_J$, the two finite sets
\[D_J^1[m]  := \big\lbrace k\in\N : 2^{J(1-a)} \leq k \leq m \big\rbrace \]
and
\begin{equation}\label{eqn:epaissi}
\mathcal{J}_J^1[m]  := \big\lbrace \mathbf{k} \in (D_J^1[m])^q \, : \, \max_{1 \leq  \ell, \ell' \leq d} |k_\ell - k_{\ell'}| \leq 2^{\varepsilon J} \big\rbrace.
\end{equation}
For all $m \in \mathcal{I}_J$, we set

\begin{equation}\label{eqn:termessimules}
 s_{m,J} := 2^{-J H} \sum_{\mathbf{k} \in \mathcal{J}_J^1[m]} \sigma_{J,\mathbf{k}}^{(q,H)} \int_\R \prod_{\ell=1}^q \Phi_\Delta^{(\delta)}(s-k_\ell) \, ds,
\end{equation}
where each random variable $\sigma_{J,\mathbf{k}}^{(q,H)}$ is defined as
\begin{equation}\label{eqn:sigma}
\sigma_{J,\mathbf{k}}^{(q,H)}:=\sum_{n=0}^{\lfloor q/2 \rfloor} (-1)^n \sum_{P \in \mathcal{P}_n^{(q)} }\prod_{r=1}^n \mathbf{E}[Z_{J,\ell_{k_r}}^{(\delta)}Z_{J,\ell_{k_r'}}^{(\delta)}] \prod_{s=n+1}^{q-n} Z_{J,\ell_{k_s''}}^{(\delta)}.
\end{equation}
One mentions that $\mathcal{P}_n^{(q)}$ in \eqref{eqn:sigma} denotes the finite set of all partitions of $\{1,\dots,q\}$ with $n$ non ordered pairs and $q-2n$ singletons and the indices $k_r$, $k_r'$ and $k_s'''$ are such that
\[ P= \left\{\{k_1,k_1'\},\dots,\{k_n,k_n'\},\{k_{n+1}''\},\dots,\{k_{q-n}''\} \right\}.\]

\begin{remark}
\label{rem:simZ}
Sample paths of the Hermite process $(Z_s, s\in I)$, on any compact interval $I \subset \R_+$, are simulated by using the restriction to $I$ of the piecewise linear continuous stochastic process
$\{\widetilde{S}_{J}^{(q,H)}(s)\}_{s \in \R_+}$ defined, for all $s \in \R_+$, as
\begin{equation}\label{simulationprocess}
\widetilde{S}_{J}^{(q,H)}(s)=\frac{s_{m_0,J}}{|\widetilde{\lambda}_{J}^{(a)}|} s\mathbbm{1}_{\widetilde{\lambda}_{J}^{(a)}}(s)+\sum_{m \in \mathcal{I}_J} \Big (2^J\left(s_{m+1,J}-s_{m,J}\right)\big (s-(m 2^{-J}+2^{-aJ})\big)+ s_{m,J} \Big) \mathbbm{1}_{\lambda_{m,J}^{(a)}}(s), 
\end{equation}
where 
\[
m_0:= \inf \mathcal{I}_J, \quad\widetilde{\lambda}_{J}^{(a)}:= [0,m_02^{-J}+2^{-aJ}), \quad |\widetilde{\lambda}_{J}^{(a)}|=m_02^{-J}+2^{-aJ},
\]
and
\[
\lambda_{m,J}^{(a)}:=[m 2^{-J}+2^{-aJ},(m+1) 2^{-J}+2^{-aJ}),\quad\mbox{for all $m \in\mathcal{I}_J $.} 
\]
Observe that 
\begin{equation}
\label{rem:simZ:e1}
\widetilde{S}_{J}^{(q,H)}\big(m 2^{-J}+2^{-aJ}\big)=s_{m,J}, \quad\mbox{for every $m \in\mathcal{I}_J $.} 
\end{equation}
The validity of this simulation procedure is guaranteed by \cite[Theorem 2.12]{AHL2}. It states that, for any compact interval $I \subset \R_+$, there exists an almost surely finite random variable $\widetilde{C}$ (depending on $I$) for which one has, almost surely, for each $J \in \N$,
\begin{equation}\label{rateofconvergence}
\| Z-\widetilde{S}_{J}^{(q,H)} \|_{I,\infty}  \leq \widetilde{C} J^{\frac{q}{2}} 2^{-J(H-\frac{1}{2})},
\end{equation}
where $\|\cdot\|_{I,\infty}$ denotes the uniform norm on $I$.
\end{remark}


\begin{remark}
\label{rem:restS}
Let $t\in\R_+$, when $s\in [0,t]$, the sum $\widetilde{S}_{J}^{(q,H)}(s)$ defined in \eqref{simulationprocess} can be rewritten as 
\begin{equation}
\label{rem:restS:e1}
\widetilde{S}_{J}^{(q,H)}(s)=\frac{s_{m_0,J}}{|\widetilde{\lambda}_{J}^{(a)}|} s\mathbbm{1}_{\widetilde{\lambda}_{J}^{(a)}}(s)+\sum_{m=m_0}^{m_t} \Big (2^J\left(s_{m+1,J}-s_{m,J}\right)\big (s-(m 2^{-J}+2^{-aJ})\big)+ s_{m,J} \Big) \mathbbm{1}_{\lambda_{m,J}^{(a)}}(s), 
\end{equation}
where 
\begin{equation}
\label{rem:restS:e2} 
m_{t}=\max \{m \in \mathcal{I}_J \, : \, m2^{-J}+2^{-aJ} \leq  t \}.
\end{equation}
Notice that when the set $\{m \in \mathcal{I}_J \, : \, m2^{-J}+2^{-aJ} \leq  t \}$ is empty then, by convention, $\sum_{m=m_0}^{m_t}\cdots=0$.
\end{remark}

\begin{definition}
\label{def:appr2}
For each fixed $J\in\N$, the piecewise linear continuous stochastic process $X_J=\{X_J (t)\}_{t\in [0,1]}$ is defined, for all $t\in [0,1]$ such that $\{m \in \mathcal{I}_J \, : \, m2^{-J}+2^{-aJ} \leq  t \}\neq \emptyset$, as
\begin{align}
\label{eqn:approxofstocint}
&X_J (t)=h(0) s_{m_0,J}+h(m_t2^{-J}+2^{-aJ})\Big (2^J\left(s_{m_t+1,J}-s_{m_t,J}\right)\big (t-(m_t 2^{-J}+2^{-aJ})\big) \Big) \nonumber \\ 
&\hspace{1.5cm}+\sum_{k=m_0}^{m_t-1}  h(k2^{-J}+2^{-aJ}) \left(s_{k+1,J}-s_{k,J}\right),
\end{align} 
where $m_t$ is as in Remark \ref{rem:restS}, and with the convention that if $m_t=m_0$ then $\sum_{k=m_0}^{m_t-1}\cdots=0$. Moreover, when $\{m \in \mathcal{I}_J \, : \, m2^{-J}+2^{-aJ} \leq  t \}= \emptyset$, then $X_J(t)$ is defined as
\begin{equation}
\label{eqn:approxofstocint-b} 
X_J (t)=h(0)\frac{s_{m_0,J}}{|\widetilde{\lambda}_{J}^{(a)}|} t .
\end{equation}
\end{definition}
The following proposition shows that when the order $\alpha\in (0,1]$ of the H\"older regularity of the function $h$ is large enough so that one has 
\[
\alpha+H>3/2,
\]
then the sequences of processes $(X_J)_{J\in\N}$ converges almost surely to the process $X$ for the uniform norm on the interval $[0,1]$, also it provides an  estimate of the rate of convergence. 
\begin{prop}
\label{prop:cvX}
Let $X$ and $X_{J}$ be given by (\ref{x:bis}) and Definition \ref{def:appr2}, respectively. For all arbitrarily small $\varepsilon>0$, there exists a positive finite random variable $C$ such that one has, almost surely, for all $J\in\N$,
\begin{equation}
\label{prop:cvX:e1}
\| X-X_J\|_{[0,1],\infty}\le C 2^{-J(\alpha+H-3/2-\varepsilon)}.
\end{equation}
\end{prop}
Our next goal is to show that Proposition \ref{prop:cvX} holds. To this end, we need to introduce the sequence of processes $(\widetilde{X}_J)_{J\in\N}$ defined as follows:
\begin{definition}
\label{def:appr}
For each fixed $J\in\N$, the continuous stochastic process $\{\widetilde{X}_J (t)\}_{t\in [0,1]}$ is defined, for all $t\in [0,1]$ such that $\{m \in \mathcal{I}_J \, : \, m2^{-J}+2^{-aJ} \leq  t \}\neq \emptyset$, as
\begin{align} 
\label{def:appr:e1}
&\widetilde{X}_J (t)=h(0)Z(m_02^{-J}+2^{-aJ})+h(m_t2^{-J}+2^{-aJ})\Big (Z(t)-Z(m_t2^{-J}+2^{-aJ}) \Big) \nonumber \\ 
&\hspace{1.5cm}+\sum_{k=m_0}^{m_t-1}  h(k2^{-J}+2^{-aJ})\Big (Z\big((k+1)2^{-J}+2^{-aJ}\big)-Z\big (k2^{-J}+2^{-aJ}\big),
\end{align} 
where $m_t$ is as in Remark \ref{rem:restS}, and with the convention that if $m_t=m_0$ then $\sum_{k=m_0}^{m_t-1}\cdots=0$. Moreover, when $\{m \in \mathcal{I}_J \, : \, m2^{-J}+2^{-aJ} \leq  t \}= \emptyset$, then $\widetilde{X}_J(t)$ is defined as
\begin{equation} 
\label{def:appr:e2}
\widetilde{X}_J (t)=h(0)Z(t).
\end{equation}
\end{definition}
{\bf Proof of Proposition \ref{prop:cvX}:} Proposition \ref{prop:cvX} is a straightforward consequence of the triangle inequality and the following two lemmas.
\qed
\begin{lemma}
\label{lem:X-XtJ}
For all arbitrarily small $\varepsilon>0$, there is a positive finite random variable $C'$ such that one has, almost surely, for all $J\in\N$,
\begin{equation}
\label{lem:X-XtJ:e1}
\| X-\widetilde{X}_J\|_{[0,1],\infty}\le C' 2^{-J(\alpha+H-1-\varepsilon)}.
\end{equation}
\end{lemma}
\begin{lemma}
\label{lem:XJ-XtJ}
For all arbitrarily small $\varepsilon>0$, there is a positive finite random variable $C''$ such that one has, almost surely, for all $J\in\N$,
\begin{equation}
\label{lem:XJ-XtJ:e1}
\| X_J-\widetilde{X}_J\|_{[0,1],\infty}\le C'' 2^{-J(\alpha+H-3/2-\varepsilon)}.
\end{equation}
\end{lemma}
{\bf Proof of Lemma \ref{lem:X-XtJ}: } In view of the assumption $(\mathcal{A})$ and the fact that sample paths of $Z$ are almost surely $(H-\varepsilon)$-H\"older continuous functions on the interval $[0,1]$,
one knows from the Young-Loeve inequality  (see Section~1.3 of \cite{LCL04}) that there is a positive finite random $C_1$, depending only $\alpha$, $H$ and $\varepsilon$, such that one has, almost surely, for all $t_1, t_2\in [0,1]$ with $t_1\le t_2$,
\begin{equation}
\label{lem:X-XtJ:e2}
\Big|\int_{t_1} ^{t_2} h(s) dZ(s)-h(t_1)\big(Z(t_2)-Z(t_1)\big)\Big| \le C_1 (t_2-t_1)^{\alpha+H-\varepsilon}.
\end{equation}
From now on $t\in [0,1]$ is arbitrary and fixed. In the sequel, we study the following three cases: first case the set $\{m \in \mathcal{I}_J \, : \, m2^{-J}+2^{-aJ} \leq  t \}$ is empty, second case it only contains $m_0$, third case it contains $m_0$ and other elements.
In the first case, using \eqref{x:bis}, \eqref{def:appr:e2}, the equality $Z(0)=0$, \eqref{lem:X-XtJ:e2}, the inequality $t<m_02^{-J}+2^{-aJ}$, the inequality $m_0\le 2^{J(1-a)}+1$ and the inequality $a(\alpha+H-\varepsilon)>\alpha+H-1-\varepsilon$, one gets, almost surely, that
\begin{equation}
\label{lem:X-XtJ:e3}
\big |X(t)-\widetilde{X}_J (t)\big|\le C_1\, t^{\alpha+H-\varepsilon}\le C_1\big (m_02^{-J}+2^{-aJ}\big)^{\alpha+H-\varepsilon}\le 9C_1 2^{-a(\alpha+H-\varepsilon)J}\le 9C_1 2^{-(\alpha+H-1-\varepsilon)J}.
\end{equation} 
In the second case, it follows from \eqref{x:bis}, \eqref{def:appr:e1}, the equality $Z(0)=0$, \eqref{lem:X-XtJ:e2}, the triangle inequality, the third and the fourth inequalities in \eqref{lem:X-XtJ:e3}, and the inequality $t\ge (m_0+1)2^{-J}+2^{-aJ}$ that
\begin{align}
\label{lem:X-XtJ:e4}
& \big |X(t)-\widetilde{X}_J (t)\big|\le\bigg | \int_{0} ^{m_02^{-J}+2^{-aJ}} h(s) dZ(s)-h(0)Z(m_02^{-J}+2^{-aJ})\bigg |\nonumber\\
&\hspace{1cm}+\bigg |\int_{m_02^{-J}+2^{-aJ}} ^{t} h(s) dZ(s)-h(m_0 2^{-J}+2^{-aJ})\Big (Z(t)-Z(m_0 2^{-J}+2^{-aJ}) \Big)\bigg | \nonumber\\
&\le C_1\big (m_02^{-J}+2^{-aJ}\big)^{\alpha+H-\varepsilon}+C_1\big (t-m_02^{-J}-2^{-aJ}\big)^{\alpha+H-\varepsilon}\nonumber\\
&\le 9C_1 2^{-a(\alpha+H-\varepsilon)J}+C_1 2^{-(\alpha+H-\varepsilon)J}\le 10C_1 2^{-(\alpha+H-1-\varepsilon)J}.
\end{align}
In the third case, one can derive from \eqref{x:bis}, \eqref{def:appr:e1}, the equality $Z(0)=0$, \eqref{lem:X-XtJ:e2}, the triangle inequality, the third and the fourth inequalities in \eqref{lem:X-XtJ:e3}, the inequality $t\ge (m_t+1)2^{-J}+2^{-aJ}$ and the inequality $m_t<2^J$ that
\begin{align}
\label{lem:X-XtJ:e5}
& \big |X(t)-\widetilde{X}_J (t)\big|\le\bigg | \int_{0} ^{m_02^{-J}+2^{-aJ}} h(s) dZ(s)-h(0)Z(m_02^{-J}+2^{-aJ})\bigg |\nonumber\\
&\hspace{1cm}+\bigg |\int_{m_t 2^{-J}+2^{-aJ}} ^{t} h(s) dZ(s)-h(m_t 2^{-J}+2^{-aJ})\Big (Z(t)-Z(m_t2^{-J}+2^{-aJ}) \Big)\bigg | \nonumber\\
&\hspace{1cm}+\sum_{k=m_0}^{m_t-1}\bigg |\int_{k2^{-J}+2^{-aJ}} ^{(k+1)2^{-J}+2^{-aJ}} h(s) dZ(s)\nonumber\\
&\hspace{3cm}-h(k 2^{-J}+2^{-aJ})\Big (Z\big ((k+1)2^{-J}+2^{-aJ}\big)-Z\big(k2^{-J}+2^{-aJ}\big) \Big)\bigg | \nonumber\\
&\le C_1\big (m_02^{-J}+2^{-aJ}\big)^{\alpha+H-\varepsilon}+C_1\big (t-m_t 2^{-J}-2^{-aJ}\big)^{\alpha+H-\varepsilon}+(m_t-m_0) C_1 2^{-(\alpha+H-\varepsilon)J}\nonumber\\
&\le 9C_1 2^{-a(\alpha+H-\varepsilon)J}+C_1 2^{-(\alpha+H-\varepsilon)J}+C_1 2^{-(\alpha+H-1-\varepsilon)J}\le 11C_1 2^{-a(\alpha+H-\varepsilon)J}.
\end{align}
Finally, setting $C'=11C_1$ and putting together \eqref{lem:X-XtJ:e3} to \eqref{lem:X-XtJ:e5}, one obtains \eqref{lem:X-XtJ:e1}.
\qed\\
\\
{\bf Proof of Lemma \ref{lem:XJ-XtJ}: } Throughout the proof $J\in\N$ and $t\in [0,1]$ are arbitrary and fixed. In the sequel, we study the following three cases: first case the set $\{m \in \mathcal{I}_J \, : \, m2^{-J}+2^{-aJ} \leq  t \}$ is empty, second case it only contains $m_0$, third case it contains $m_0$ and other elements.
In the first case, using \eqref{def:appr:e2}, \eqref{eqn:approxofstocint-b}, Remark \ref{rem:restS} and \eqref{rateofconvergence}, one obtains that 
\begin{align}
\label{lem:XJ-XtJ:e2}
\big | X_J (t)-\widetilde{X}_J(t)\big| =|h(0)| \big | \widetilde{S}_{J}^{(q,H)}(t)-Z(t)\big|\le \|h\|_{[0,1],\infty}\widetilde{C} J^{\frac{q}{2}} 2^{-J(H-\frac{1}{2})}.
\end{align} 
In the second case, one can derive from \eqref{def:appr:e1}, \eqref{eqn:approxofstocint}, the triangle inequality, \eqref{rem:simZ:e1}, \eqref{rem:restS:e1} and \eqref{rateofconvergence} that 
\begin{align}
\label{lem:XJ-XtJ:e3}
&\big | X_J (t)-\widetilde{X}_J(t)\big| \le |h(0)| \big | \widetilde{S}_{J}^{(q,H)}\big (m_02^{-J}+2^{-aJ}\big)-Z\big (m_02^{-J}+2^{-aJ}\big)\big|\nonumber\\
&\hspace{3mm}+\big |h\big(m_02^{-J}+2^{-aJ}\big)\big|\bigg|\Big (\widetilde{S}_{J}^{(q,H)}(t)-\widetilde{S}_{J}^{(q,H)}\big (m_02^{-J}+2^{-aJ}\big)\Big)-\Big (Z(t)-Z\big (m_02^{-J}+2^{-aJ}\big)\Big)\bigg| \nonumber\\
&\le 2\|h\|_{[0,1],\infty} \big | \widetilde{S}_{J}^{(q,H)}\big (m_02^{-J}+2^{-aJ}\big)-Z\big (m_02^{-J}+2^{-aJ}\big)\big|+\|h\|_{[0,1],\infty} \big |\widetilde{S}_{J}^{(q,H)}(t)-Z(t)\big|\nonumber\\
&\le 3\|h\|_{[0,1],\infty}\widetilde{C} J^{\frac{q}{2}} 2^{-J(H-\frac{1}{2})}. 
\end{align}
In the third case, it follows from \eqref{def:appr:e1}, \eqref{eqn:approxofstocint}, the triangle inequality, \eqref{rem:simZ:e1}, \eqref{rem:restS:e1} and \eqref{rateofconvergence} that 
\begin{align}
\label{lem:XJ-XtJ:e4}
& \big | X_J (t)-\widetilde{X}_J(t)\big| \le |h(0)| \big | \widetilde{S}_{J}^{(q,H)}\big (m_02^{-J}+2^{-aJ}\big)-Z\big (m_02^{-J}+2^{-aJ}\big)\big|\nonumber\\
&\hspace{3mm}+\big |h\big(m_t 2^{-J}+2^{-aJ}\big)\big|\bigg|\Big (\widetilde{S}_{J}^{(q,H)}(t)-\widetilde{S}_{J}^{(q,H)}\big (m_t 2^{-J}+2^{-aJ}\big)\Big)-\Big (Z(t)-Z\big (m_t 2^{-J}+2^{-aJ}\big)\Big)\bigg| \nonumber\\
&\hspace{3mm}+\Big |\sum_{k=m_0}^{m_t-1}h\big(k 2^{-J}+2^{-aJ}\big)\big (\delta_{J,k+1}-\delta_{J,k}\big)\Big|\nonumber\\
&\le 3\|h\|_{[0,1],\infty}\widetilde{C} J^{\frac{q}{2}} 2^{-J(H-\frac{1}{2})} +\Big |\sum_{k=m_0}^{m_t-1}h\big(k 2^{-J}+2^{-aJ}\big)\big (\delta_{J,k+1}-\delta_{J,k}\big)\Big|,
\end{align}
where 
\begin{equation}
\label{lem:XJ-XtJ:e5}
\delta_{J,k}=\widetilde{S}_{J}^{(q,H)}\big(k 2^{-J}+2^{-aJ}\big)-Z\big(k 2^{-J}+2^{-aJ}\big),\quad\mbox{for all $k\in\{m_0,\ldots,m_t\}$.}
\end{equation}
Moreover, one has that 
\begin{align*}
& \sum_{k=m_0}^{m_t-1}h\big(k 2^{-J}+2^{-aJ}\big)\big (\delta_{J,k+1}-\delta_{J,k}\big)\\
&=\sum_{k=m_0+1}^{m_t}h\big((k-1)2^{-J}+2^{-aJ}\big)\delta_{J,k}-\sum_{k=m_0}^{m_t-1}h\big(k 2^{-J}+2^{-aJ}\big)\delta_{J,k}\\
&=h\big((m_t-1)2^{-J}+2^{-aJ}\big)\delta_{J,m_t}-h\big(m_0 2^{-J}+2^{-aJ}\big)\delta_{J,m_0}\\
&\hspace{1cm}-\sum_{k=m_0+1}^{m_t-1}\Big( h\big(k 2^{-J}+2^{-aJ}\big)-h\big((k-1)2^{-J}+2^{-aJ}\big)\Big)\delta_{J,k}.
\end{align*}
Thus, the triangle inequality, \eqref{lem:XJ-XtJ:e5}, \eqref{rateofconvergence}, the $\alpha$-H\"older continuity of $h$ and the inequality $m_t<2^J$ imply that
\begin{align}
\label{lem:XJ-XtJ:e6}
&\Big |\sum_{k=m_0}^{m_t-1}h\big(k 2^{-J}+2^{-aJ}\big)\big (\delta_{J,k+1}-\delta_{J,k}\big)\Big|\le 2\|h\|_{[0,1],\infty}\widetilde{C} J^{\frac{q}{2}} 2^{-J(H-\frac{1}{2})}\\
&\hspace{3mm}+\widetilde{C} J^{\frac{q}{2}} 2^{-J(H-\frac{1}{2})} \sum_{k=m_0+1}^{m_t-1}\Big |h\big(k 2^{-J}+2^{-aJ}\big)-h\big((k-1)2^{-J}+2^{-aJ}\big)\Big|\nonumber
&\le C_1 2^{-J(\alpha+H-3/2-\varepsilon)},
\end{align} 
where $C_1$ is a positive finite random variable not depending on $J$ and $t$.
Finally, putting together \eqref{lem:XJ-XtJ:e2}, \eqref{lem:XJ-XtJ:e3}, \eqref{lem:XJ-XtJ:e4} and \eqref{lem:XJ-XtJ:e6}, one obtains \eqref{lem:XJ-XtJ:e1}
\qed

In this work, we conducted experiments for Hermite processes of order $q = 1, 2, 3$. Note that, in these cases, the expression \eqref{eqn:sigma} take the more manageable form
\[ \sigma_{J,\mathbf{k}}^{(q,H)} := \begin{cases}
Z_{J,k}^{(\delta)} & \text{if } q=1 \\
Z_{J,k_{1}}^{(\delta)} Z_{J,k_{2}}^{(\delta)} -\mathbb{E}[Z_{J,k_{1}}^{(\delta)} Z_{J,k_{2}}^{(\delta)} ] & \text{if } q=2 \\
Z_{J,k_{1}}^{(\delta)}Z_{J,k_{2}}^{(\delta)}Z_{J,k_{3}}^{(\delta)}-\mathbb{E}[Z_{J,k_{1}}^{(\delta)}Z_{J,k_{2}}^{(\delta)}]Z_{J,k_{3}}^{(\delta)} & \text{if } q=3 \\
-\mathbb{E}[Z_{J,k_{1}}^{(\delta)}Z_{J,k_{3}}^{(\delta)}]Z_{J,k_{2}}^{(\delta)}-\mathbb{E}[Z_{J,k_{2}}^{(\delta)}Z_{J,k_{3}}^{(\delta)}]Z_{J,k_{1}}^{(\delta)} & 
\end{cases}\]
and explicit algorithms of simulations for these processes are given in \cite{AHL2}.  Note that this last paper and \cite{LT} also present numerical experiments for various estimators for the Hurst parameter of Hermite processes. In particular, the estimator based on a modified quadratic variation defined in \cite{AT} and the estimator based on a modified wavelet variation defined in \cite{LT} were considered. Both estimators are build using an approach similar to the one explored in the current paper. In the context of Hurst parameter estimation, it was observed that the higher the resolution, the better the estimator performs. For this reason, the numerical experiments we carry out below are performed at scale $J=18$, which represents a good compromise between approximation quality and computation time for generating a large sample of simulations. Therefore, in order not to work with an estimator that requires considering more data than those generated by a single simulation, we propose to use $N=17$. 

The experiments presented various estimations obtained for the Fractional Brownian motion (FBM), Rosenblatt process (RP) and Hermite process of order 3 (HP) simulated with prescribed Hurst parameters: $0.6$, $0.7$, $0.8$ and $0.9$. In each cases, we simulated 100 trajectories with  $a = 0.99$ and $\varepsilon=10^{-3}$. The notations $m$ and $s$ stand respectively for the mean of the estimations and their standard deviation. Moments of order $p \neq 2$ for the Hermite distribution of order $q>1$ are not explicitly known. For this reason, we only consider the case $p=2$ in the Tables below. 

First experiments are carried out with the identity function $h(s)=s$ for which the expected value of the limit in \eqref{eqn:convergencelimite} is $\frac{1}{3}$. The estimations of the integrated volatility obtained when we choose $\gamma=0.8$ are presented in Table 1 while, in Table 2, we use $\gamma=0.95$. We observe that, in general, the estimation is more accurate when the Hurst index $H$ is large. This can certainly be explained by equation \eqref{prop:cvX:e1}, which entails that the larger $H$ is, the more precise the approximation of the stochastic integral  \eqref{x} by the process \eqref{eqn:approxofstocint} becomes. We also note that, on average, the estimators defined with $\gamma = 0.8$ and $\gamma = 0.95$ have similar performance, but the standard deviation of the estimates obtained with $\gamma = 0.95$ is smaller than that with $\gamma = 0.8$. This can perhaps be explained in light of definition \eqref{snpx}. The larger $\gamma$ is, the more terms there are in this sum, and therefore the more the estimator makes use of the data available to it. For this reason, we choose $\gamma = 0.95$ for the remainder of our experiments.

In Table 3, we present the estimations for the integrated volatility in the case $h(s)=s^3$, where the expected value is $\frac{1}{7}\approx 0.142$, and Table 4 concerns the case $h(s)=e^s$ with expected value $\frac{e^2-1}{2}\approx3.194$. We again see that the estimator performs quite well. The rate of convergence obtained in (\ref{rateofconvergence}) certainly explains why the performance of the estimation appears to be better when $q$ is smaller.

Proposition \ref{prop:cvX} states that as soon as $\alpha+H>3/2$, the sequence of stochastic process $(\{X_J(t)\}_{t \in [0,1]})_J$ almost surely uniformly converges to the process $\{X_t\}_{t \in [0,1]}$ defined by \eqref{x}. It validates the approach to simulate the data using a process $\{X_J(t)\}_{t \in [0,1]}$ as soon as $\alpha+H>3/2$. In order to investigate the sharpness of this assumption, we estimate the integrated volatility in the case $h(x)=\sqrt{x}$. The accuracy of the estimations in Table 5 seems to show that our simulation approach remains valid even if the assumption $\alpha+H>3/2$ does not hold. A formal proof of this statement would require a more in-depth analysis of the convergence in Proposition \ref{prop:cvX}. A related interesting research project would be to investigate whether other relevant simulation procedures could be developed for the stochastic integral of the form \eqref{x}. A promising approach would be to extend the wavelet-type expansion obtained in \cite{AHL1} in this context and exploit it in the same spirit as in \cite{AHL2}.

\begin{table}[h]
\begin{center}
\begin{tabular}{|c|c|c|c|c|c|}
\hline 
\multicolumn{2}{|c|}{ }& $0.6$ & $0.7$ & $0.8$ & $0.9$   \\ 
\hline 
\multirow{2}{*}{1} & m &0.295  & 0.291  & 0.298 &  0.307 \\
\cline{2-6} 
 & s & 0.024 & 0.025 & 0.024 & 0.023  \\ 
\hline 
\multirow{2}{*}{2} & m &0.374 & 0.370 & 0.331 & 0.344\\ 
\cline{2-6}
 & s & 0.065& 0.075 &0.071&  0.070 \\ 
\hline
\multirow{2}{*}{3} & m & 0.392 &0.388& 0.347 &  0.343 \\ 
\cline{2-6}
 & s &0.159 & 0.212 &0.161& 0.053 \\ 
\hline 
\end{tabular} 
\caption{Integrated volatility estimation for $h(s)=s$ with $J=18$, $N=17$ and $\gamma=0.8$}
\end{center}
\end{table}

\begin{table}[h]
\begin{center}
\begin{tabular}{|c|c|c|c|c|c|}
\hline 
\multicolumn{2}{|c|}{ }& $0.6$ & $0.7$ & $0.8$ & $0.9$   \\ 
\hline 
\multirow{2}{*}{1} & m &0.291  & 0.290  & 0.296 &  0.306 \\
\cline{2-6} 
 & s & 0.003 & 0.003 & 0.005 & 0.013  \\ 
\hline 
\multirow{2}{*}{2} & m &0.386& 0.364 & 0.345 & 0.351\\ 
\cline{2-6}
 & s & 0.017& 0.018 &0.026&  0.034 \\ 
\hline
\multirow{2}{*}{3} & m & 0.412 &0.383& 0.350 &  0.346 \\ 
\cline{2-6}
 & s &0.049 & 0.051 &0.046& 0.026 \\ 
\hline 
\end{tabular} 
\caption{Integrated volatility estimation for $h(s)=s$ with $J=18$, $N=17$ and $\gamma=0.95$}
\end{center}
\end{table}

\begin{table}[h]
\begin{center}
\begin{tabular}{|c|c|c|c|c|c|}
\hline 
\multicolumn{2}{|c|}{ }& $0.6$ & $0.7$ & $0.8$ & $0.9$   \\ 
\hline 
\multirow{2}{*}{1} & m &0.125  & 0.124   &  0.127 &  0.131 \\
\cline{2-6} 
 & s & 0.002 & 0.002 & 0.003& 0.008\\ 
\hline 
\multirow{2}{*}{2} & m &0.166 & 0.155 & 0.147 & 0.147\\ 
\cline{2-6}
 & s & 0.011& 0.012 &0.019 &  0.021 \\ 
\hline
\multirow{2}{*}{3} & m & 0.174&0.164& 0.148 & 0.147 \\ 
\cline{2-6}
 & s &0.032& 0.033 & 0.031 & 0.017 \\ 
\hline 
\end{tabular} 
\caption{Integrated volatility estimation  for $h(s)=s^3$ with $J=18$, $N=17$ and $\gamma=0.95$}
\end{center}
\end{table}

\begin{table}[h]
\begin{center}
\begin{tabular}{|c|c|c|c|c|c|}
\hline 
\multicolumn{2}{|c|}{ }& $0.6$ & $0.7$ & $0.8$ & $0.9$   \\ 
\hline 
\multirow{2}{*}{1} & m &2.983 & 2.803  & 2.926 &  3.268  \\
\cline{2-6} 
 & s & 0.032 & 0.027 & 0.038 & 0.093\\ 
\hline 
\multirow{2}{*}{2} & m &3.696 &3.501 & 3.323 & 3.397\\ 
\cline{2-6}
 & s & 0.135& 0.140 &0.172 &  0.242 \\ 
\hline
\multirow{2}{*}{3} & m & 3.974 &3.682& 3.394 & 3.335 \\ 
\cline{2-6}
 & s &0.355& 0.346 & 0.317 & 0.173 \\ 
\hline 
\end{tabular} 
\caption{Integrated volatility estimation  for $h(s)=\exp(s)$ with $J=18$, $N=17$ and $\gamma=0.95$}
\end{center}
\end{table}

\begin{table}[h]
\begin{center}
\begin{tabular}{|c|c|c|c|c|c|}
\hline 
\multicolumn{2}{|c|}{ }& $0.6$ & $0.7$ & $0.8$ & $0.9$   \\ 
\hline 
\multirow{2}{*}{1} & m &0.437  & 0.436  & 0.472 &  0.459  \\
\cline{2-6} 
 & s & 0.005& 0.004 & 0.006 & 0.015\\ 
\hline 
\multirow{2}{*}{2} & m &0.578 & 0.547& 0.519 & 0.531\\ 
\cline{2-6}
 & s & 0.020 & 0.023 &0.027 &  0.040 \\ 
 \hline
\multirow{2}{*}{3} & m & 0.622 &0.575& 0.530 & 0.521 \\ 
\cline{2-6}
 & s &0.057& 0.058 & 0.052& 0.029 \\ 
\hline 
\end{tabular} 
\caption{Integrated volatility estimation  for $h(s)=\sqrt{s}$ with $J=18$, $N=17$ and $\gamma=0.95$}
\end{center}
\end{table}

\clearpage

\section*{Acknowledgements}
A. Ayache acknowledges partial support from the CDP C2EMPI, as well as the French State under the France-2030 programme, the University of Lille, the Initiative of Excellence of the University of Lille and the European Metropolis of Lille by their funding and support of the R-CDP-24-004-C2EMPI project. Also, he has been partially supported by the Australian Research Council's Discovery Projects funding scheme (project number  DP220101680).

L. Loosveldt acknowledges support from the FRNS Research Credit grant J.0136.26 ``NOISE''.

 C. Tudor  acknowledges support from   the  ANR project SDAIM 22-CE40-0015, the MATHAMSUD grant 24-MATH-04 SDE-EXPLORE,  and  the Ministry of Research, Innovation and Digitalization (Romania), grant CF-194-PNRR-III-C9-2023.

	\section{Appendix: Multiple stochastic integrals and the Malliavin derivative}\label{app}

The basic tools from the analysis on Wiener space are presented in this section. We will focus on some elementary facts about multiple stochastic integrals. We refer to \cite{N} or \cite{NP-book} for a complete review on the topic. 

Consider ${\mathcal{H}}$ a real separable infinite-dimensional Hilbert space
with its associated inner product ${\langle
	\cdot,\cdot\rangle}_{\mathcal{H}}$, and $(B (\varphi),
\varphi\in{\mathcal{H}})$ an isonormal Gaussian process on a
probability space $(\Omega, {\mathfrak{F}}, \mathbb{P})$, which is a
centred Gaussian family of random variables such that
$\mathbf{E}\left( B(\varphi) B(\psi) \right) = {\langle\varphi,
	\psi\rangle}_{{\mathcal{H}}}$ for every
$\varphi,\psi\in{\mathcal{H}}$. Denote by $I_{q}$ ($q\geq 1$) the $q$th multiple
stochastic integral with respect to $B$, which is an
isometry between the Hilbert space ${\mathcal{H}}^{\odot q}$
(symmetric tensor product) equipped with the scaled norm
$\sqrt{q!}\,\Vert\cdot\Vert_{{\mathcal{H}}^{\otimes q}}$ and
the Wiener chaos of order $q$, which is defined as the closed linear
span of the random variables $H_{q}(B(\varphi))$ where
$\varphi\in{\mathcal{H}},\;\Vert\varphi\Vert_{{\mathcal{H}}}=1$ and
$H_{q}$ is the Hermite polynomial of degree $q\geq 1$ defined
by:
\begin{equation}\label{Hermite-poly}
	H_{q}(x)=\frac{(-1)^{q}}{q!} \exp \left( \frac{x^{2}}{2} \right) \frac{{\mathrm{d}}^{q}%
	}{{\mathrm{d}x}^{q}}\left( \exp \left(
	-\frac{x^{2}}{2}\right)\right),\;x\in \mathbb{R}.
\end{equation}
For $q=0$, 
\begin{equation}\label{27a-1}
	\mathcal{H}_{0}=\mathbb{R} \mbox{ and }I_{0}(x)=x \mbox{ for every }x\in \mathbb{R}.
\end{equation}
The isometry property of multiple integrals can be written as follows : for $p,\;q\geq
0$,\;$f\in{{\mathcal{H}}^{\otimes p}}$ and
$g\in{{\mathcal{H}}^{\otimes q}}$
\begin{equation} \mathbf{E}\Big(I_{p}(f) I_{q}(g) \Big)= \left\{
	\begin{array}{rcl}\label{iso}
		q! \langle \tilde{f},\tilde{g}
		\rangle _{{\mathcal{H}}^{\otimes q}}&&\mbox{if}\;p=q,\\
		\noalign{\vskip 2mm} 0 \quad\quad&&\mbox{otherwise,}
	\end{array}\right.
\end{equation}
where $\tilde{f}$ stands for the symmetrization of $f$. When $\mathcal{H}= L^{2}(T)$, with $T$ being an interval of $\mathbb{R}$, we have the following product formula: for $p,\;q\geq
0$,\;  $f\in{{\mathcal{H}}^{\odot p}}$ and
$g\in{{\mathcal{H}}^{\odot q}}$,  

\begin{eqnarray}\label{prod}
	I_{p}(f) I_{q}(g)&=& \sum_{r=0}^{p \wedge q} r! \binom{q}{r} \binom{p}{r} I_{p+q-2r}\left(f\otimes_{r}g\right),
\end{eqnarray}

where, for $r=0, ..., p\wedge q$, the contraction $f\otimes _{r} g$ is the function in $L ^{2}( T ^{p+q-2r}) $ given by 

\begin{equation}\label{contra}
	(f\otimes _{r} g) (t_{1},..., t_{p+q-2r})= \int_{T^{r}} f(u_{1},..., u_{r}, t_{1},..., t_{p-r}) g(u_{1},..., u_{r}, t_{p-r+1},...,t_{p+q-2r}) du_{1}...du_{r}.
\end{equation}

An useful  property of  finite sums of multiple stochastic integrals is the hypercontractivity. Namely,  for every fixed real number $p\geq 2$, there exists a universal deterministic finite constant $C_p$, such that, for any random variable $F$ of the form $F= \sum_{k=0} ^{n} I_{k}(f_{k}) $ with $f_{k}\in \mathcal{H} ^{\otimes k}$, the following inequality holds:
\begin{equation}
	\label{hyper}
	\mathbf{E}\vert F \vert ^{p} \leq C_{p} \left( \mathbf{E}F ^{2} \right) ^{\frac{p}{2}}.
\end{equation}

We denote by $D$ the Malliavin derivative operator that acts on
cylindrical random variables of the form $F=g(B(\varphi
_{1}),\ldots,B(\varphi_{n}))$, where $n\geq 1$,
$g:\mathbb{R}^n\rightarrow\mathbb{R}$ is a smooth function with
compact support and $\varphi_{i} \in {{\mathcal{H}}}$, in the following way:
\begin{equation*}
	DF=\sum_{i=1}^{n}\frac{\partial g}{\partial x_{i}}(B(\varphi _{1}),
	\ldots , B(\varphi_{n}))\varphi_{i}.
\end{equation*}
The operator $D$ is closable and it can be extended to $\mathbb{D} ^{1, 2}$ which denotes the closure of the set of cylindrical random variables with respect to the norm $\Vert \cdot\Vert _{1,2}$ defined as
\begin{equation*}
	\Vert F\Vert _{1,2} ^{2}:= \mathbf{E}\vert F\vert ^{2}+ \mathbf{E} \Vert  DF\Vert _{\mathcal{H}} ^{2}. 
\end{equation*} 
If $F=I_{p}(f)$, where $f\in \mathcal{H} ^{\odot p}$ with $\mathcal{H}= L^{2}(T)$ and $p\geq 1$, then 
$$D_{\ast}F=pI_{p-1} \left( f(\cdot, \ast)\right),$$
where $"\cdot "$ stands for $p-1$ variables.

The pseudo inverse $ (-L) ^{-1}$ of the Ornstein-Uhlenbeck operator $L$ is defined, for $F=I_{p}(f)$ with $f\in \mathcal{H} ^{\odot p}$ and $p\geq 1$, by 
\begin{equation*}
	(-L) ^{-1}F= \frac{1}{p} I_{p} (f).
\end{equation*}

At last notice that in our work, we have $\mathcal{H}= L ^{2} (\mathbb{R})$ while the role of the isonormal process $(B(\varphi), \varphi \in \mathcal{H}) $ is played by  the usual Wiener integral on $L ^{2} (\mathbb{R})$ associated with the Brownian motion $(B(y), y \in \mathbb{R})$.

\end{document}